\documentclass[
aps,%
pre,%
a4paper,%
final,%
reprint,%
showpacs,%
superscriptaddress,%
groupedaddress,%
nofootinbib
]{revtex4-1}

\usepackage{amsmath}
\usepackage{mathrsfs}  
\usepackage{amssymb}
\usepackage{amsthm}
\usepackage{bm}
\usepackage{layouts}

\usepackage[pdftex]{graphicx}
\graphicspath{{./figs/}}
\usepackage[pdftex,colorlinks=true,bookmarks=true,citecolor=blue,urlcolor=blue]{hyperref} 
\usepackage[caption=false]{subfig}

\newcommand{\field}[1]{\mathbb{#1}}
\newcommand{\fs}[1]{\mathsf{#1}}

\DeclareMathOperator*{\esssup}{ess\,sup}

\DeclareMathOperator{\diag}{diag}

\newcommand{\tp}{\intercal}

\newcommand{\ovl}[1]{\overline{#1}}

\newcommand{\bigO}[1]{\mathop{\mathscr{O}}\left(#1\right)}

\let\Re\relax
\DeclareMathOperator{\Re}{Re}
\let\Im\relax
\DeclareMathOperator{\Im}{Im}
\newcommand{\vv}[1]{\boldsymbol{#1}}
\newcommand{\vs}[1]{\boldsymbol{#1}}

\newcommand{\OP}[1]{\mathscr{#1}}

\DeclareMathOperator{\fourier}{\mathscr{F}}

\DeclareMathOperator{\sech}{sech}

\DeclareMathOperator{\sinc}{sinc}

\newcommand{\wtilde}[1]{\widetilde{#1}}

\newcommand{\et}{\textit{et~al.}}

\usepackage{enumitem}

\begin{document}
\title{Fast Inverse Nonlinear Fourier Transform}

\author{V.~Vaibhav}
\email{vishal.vaibhav@gmail.com}

\date{\today}

\begin{abstract}
This paper considers the non-Hermitian Zakharov-Shabat (ZS) scattering problem
which forms the basis for defining the SU$(2)$-nonlinear Fourier transform (NFT). The 
theoretical underpinnings of this generalization of the conventional Fourier 
transform is quite well established in the Ablowitz-Kaup-Newell-Segur (AKNS)
formalism; however, efficient numerical algorithms that could be employed in
practical applications are still unavailable. In this paper, we present two 
fast inverse NFT algorithms with $O(KN+N\log^2N)$ complexity and a 
convergence rate of $O(N^{-2})$ where $N$ is the number of samples of 
the signal and $K$ is the number of eigenvalues. These algorithms are 
realized using a new fast layer-peeling (LP) scheme ($O(N\log^2N)$) together 
with a new fast Darboux transformation (FDT) algorithm ($O(KN+N\log^2N)$) previously 
developed by the author~[Phys. Rev. E 96, 063302 (2017)]. The proposed fast inverse NFT algorithm 
proceeds in two steps: The first
step involves computing the radiative part of the potential using the fast LP
scheme for which the input is synthesized under the assumption
that the radiative potential is nonlinearly bandlimited, i.e., the continuous
spectrum has a compact support and the discrete spectrum is empty. The second 
step involves addition of bound states using the FDT algorithm. Finally, the 
performance of these algorithms is demonstrated through exhaustive numerical tests.
\end{abstract}

\pacs{%
02.30.Zz,
02.30.Ik,
42.81.Dp,
03.65.Nk
}

\maketitle

\section*{Notations}
\label{sec:notations}
The set of non-zero positive real numbers ($\field{R}$) is denoted by
$\field{R}_+$. Non-zero positive (negative) integers are denoted by
$\field{Z}_+$ ($\field{Z}_-$). For any complex number $\zeta$, 
$\Re(\zeta)$ and $\Im(\zeta)$ refer to the real
and the imaginary parts of $\zeta$, respectively. Its complex conjugate is 
denoted by $\zeta^*$. The upper-half (lower-half) of complex plane, $\field{C}$, 
is denoted by $\field{C}_+$ ($\field{C}_-$).

\section{Introduction}
The nonlinear Fourier (NF) spectrum offers a novel way of encoding information 
in optical pulses where the nonlinear effects are adequately taken into 
account as opposed to being treated as a source of distortion. This idea 
has its origin in the work of Hasegawa and Nyu~\cite{HN1993} who were the first
to propose the use of discrete eigenvalues of the NF spectrum for encoding 
information. Recent 
advances in coherent optical communication have made it possible to reconsider 
this old idea with some extensions and improvements. Extension of this scheme
consists in using additional degrees of freedom offered by the NF 
spectrum such as the norming constants and the continuous spectrum. For an overview 
of the recent progress in theoretical as well as experimental aspects of various 
optical communication methodologies that are based on the nonlinear Fourier 
transform (NFT), we refer the reader to the review article~\cite{TPLWFK2017} and
the references therein.

\begin{figure*}[!ht]
\centering
\includegraphics[scale=1]{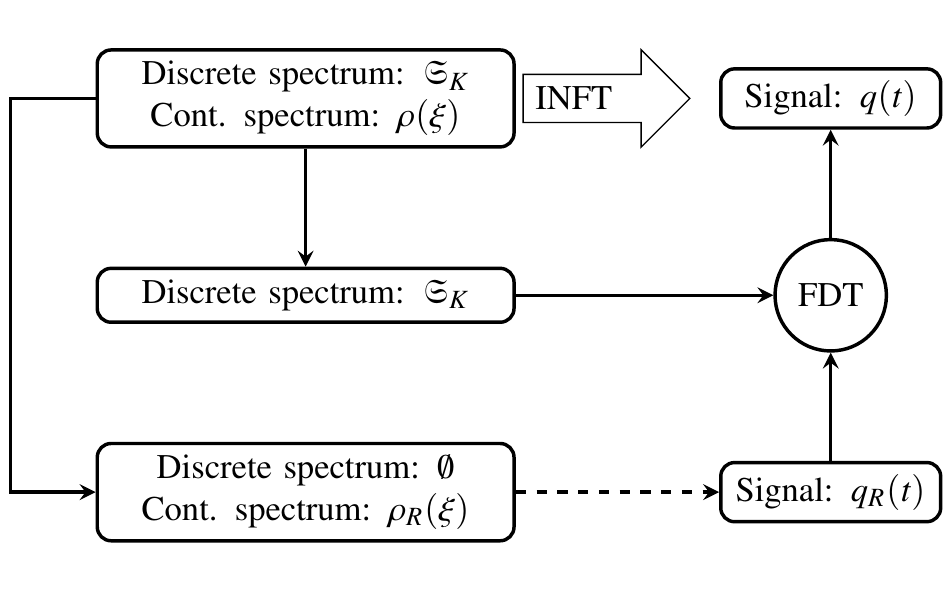}
\caption{\label{fig:schema-inft}The figure shows a schematic of the
fast inverse NFT (INFT) algorithm where the dashed line depicts
the missing part of the algorithm to be discussed in this article 
(the FDT algorithm has been reported in~\cite{V2017INFT1}). 
Here, $q_R(t)$ refers to the ``radiative'' part of the signal
$q(t)$ which is obtained as a result of removing the bound states. Note that
$q_R(t)$ has the reflection coefficient $\rho_R(\xi)$ (see
Sec.~\ref{sec:akns-sys} for the connection between $\rho_R(\xi)$ and $\rho(\xi)$).}
\end{figure*}

In order to realize any NFT-based modulation methodology, it is 
imperative to have a suitable low-complexity
NFT algorithm which forms the primary motivation behind this work. The central
idea is to use a fast version of the well-known \emph{layer-peeling} (LP) algorithm 
within the framework of an appropriate discretization scheme applied to the
Zakharov-Shabat (ZS) problem. This approach has been characterized as the
differential approach by Bruckstein~\et~\cite{BLK1985,BK1987} where fast realizations
of the LP algorithm which achieves a complexity of $\bigO{N\log^2N}$ 
for $N$ samples of the reflection data are also discussed. However, the earliest
work on fast LP is that of McClary~\cite{McClary1983} which appeared
in the geophysics literature. More recently, this method has been adopted by 
Brenne and Skaar~\cite{BS2003} in the design of grating-assisted codirectional 
couplers. However, this paper reports a complexity of $\bigO{N^2}$\footnote{In
this paper, we do not consider the method of discretization presented 
in~\cite{BS2003}; however, let us briefly mention that on account of the
piecewise constant assumption used in this work for the scattering potential, the 
order of convergence gets artificially restricted to $\bigO{N^{-1}}$. This 
problem has been remedied in~\cite{V2017INFT1} where this 
discretization scheme is termed as the \emph{split-Magnus} method.}. It is
interesting to note that, at the heart of it, all of the aforementioned versions of
LP are similar; however, the manner in which the discrete system is obtained
seem to vary. In this work, we consider the discrete system obtained as a result of
applying (exponential) trapezoidal rule to the ZS problem as discussed 
in~\cite{V2017INFT1}. 
 
The next important idea is to recognize that the 
\emph{Darboux transformation} (DT) provides a 
promising route to the most general inverse NFT 
algorithm. A fast version of DT (referred to as FDT) is developed
in~\cite{V2017INFT1} which is based on the pioneering work of 
Lubich on convolution quadrature~\cite{Lubich1994} and a 
fast LP algorithm. The schematic of the fast inverse NFT is shown 
in Fig.~\ref{fig:schema-inft} where we note that FDT is
capable of taking a \emph{seed} potential $q_R(t)$ and augmenting it by
introducing the bound states corresponding to $\mathfrak{S}_K$ (the discrete
spectrum to be introduced in Sec.~\ref{sec:akns-sys} and $K$ is the number of
bound states or eigenvalues). If $q_R(t)$ is the
\emph{radiative} part of $q(t)$, i.e., it is generated from NF spectrum which
has an empty discrete spectrum and $\rho_R(\xi)$ as the reflection coefficient, then $q(t)$
is the full inverse of the NF spectrum characterized by $\mathfrak{S}_K$ and
$\rho(\xi)$. The preliminary results of this approach 
were reported in~\cite{VW2017OFC}. In this paper, we describe two fast inverse NFT
algorithms that exhibit a complexity of $\bigO{N(K+\log^2N)}$ and a rate of 
convergence of $\bigO{N^{-2}}$ where $N$ is
the number of samples and $K$ is the number of eigenvalues (or bound states).
 
Finally, we note that the LP algorithm (irrespective of the underlying
discrete system) has the reputation of being ill-conditioned or unstable in the
presence of noise~\cite{BKK1986,SF2002} in the reflection coefficient. For optical
communication, this observation is important but not critical as the reflection
coefficient is known exactly at the stage of encoding of information at the
transmitter end. A more relevant question here, therefore, is the stability of the 
algorithm in the presence of round-off errors. We provide exhaustive numerical tests 
in order to understand the ill-conditioning effects; however, no
theoretical results for stability are provided.

This paper is organized as follows: Sec.~\ref{sec:akns-sys} discusses the basic
theory of scattering. Sec.~\ref{sec:discrete-system} introduces the discrete
framework for forward/inverse scattering, which admits of the layer-peeling
property. This section also introduces a recipe for
computing a class of signals dubbed as the \emph{nonlinearly bandlimited}
signals. Finally, the inverse NFT 
is described in Sec.~\ref{sec:fast-inverse-NFT} and the numerical results are 
presented in Sec.~\ref{sec:num-res}. Sec.~\ref{sec:final} concludes this paper.

\begin{figure*}[!ht]
\centering
\includegraphics[scale=1]{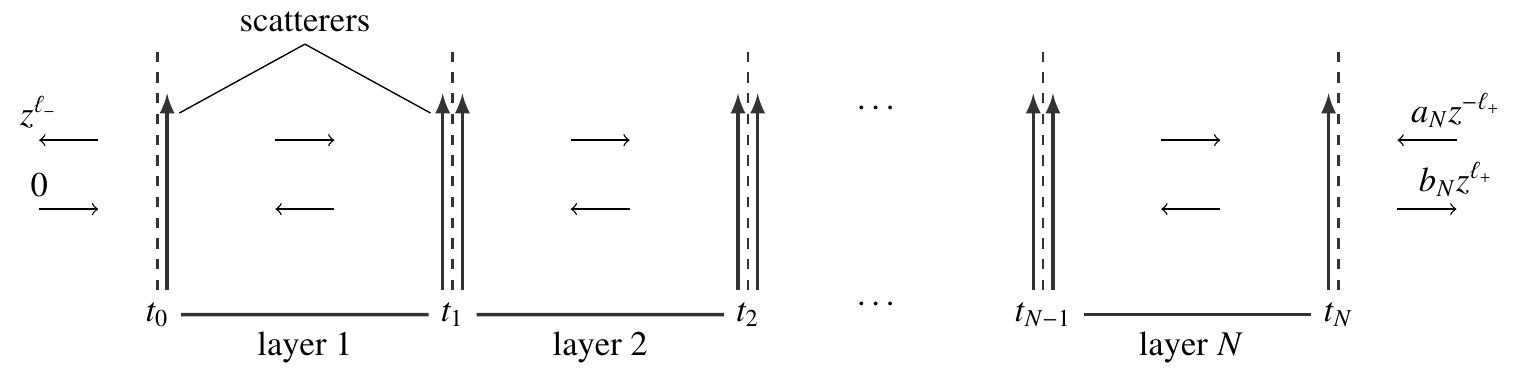}
\caption{The figure depicts the equivalent \emph{layered-media} for the discrete
scattering problem in Sec.~\ref{sec:discrete-system}. In each of the layers, the ZS-problem 
is approximated by two instantaneous scatterers and a ``free-space'' propagation 
between them.\label{fig:sm}}
\end{figure*}

\section{The AKNS System}
\label{sec:akns-sys}
The NFT of any complex-valued signal $q(t)$ is introduced 
via the associated Zakharov-Shabat scattering problem~\cite{ZS1972} 
which can be stated as follows:
Let $\zeta\in\field{R}$ and $\vv{v}=(v_1,v_2)^{\tp}\in\field{C}^2$, then 
\begin{equation}\label{eq:zs-prob}
\vv{v}_t = -i\zeta\sigma_3\vv{v}+U\vv{v}, 
\end{equation}
where $\sigma_3=\diag(1,-1)$, and, the matrix elements of $U$ are $U_{11}=U_{22}=0$ and
$U_{12}=q(t)=-U_{21}^*=-r^*(t)$. Here, $q(t)$ is identified as the \emph{scattering potential}. The 
solution of the scattering 
problem~\eqref{eq:zs-prob}, henceforth referred
to as the ZS problem, consists in finding the so called 
\emph{scattering coefficients} which are defined through 
special solutions of~\eqref{eq:zs-prob}
known as the \emph{Jost solutions}. The Jost solutions of the \emph{first kind}, denoted
by $\vs{\psi}(t;\zeta)$, has the asymptotic behavior 
$\vs{\psi}(t;\zeta)e^{-i\zeta t}\rightarrow(0,1)^{\tp}$ as $t\rightarrow\infty$. 
The Jost solutions of the \emph{second kind}, denoted by $\vs{\phi}(t,\zeta)$, has 
the asymptotic behavior $\vs{\phi}(t;\zeta)e^{i\zeta t}\rightarrow(1,0)^{\tp}$ 
as $t\rightarrow-\infty$. The so-called scattering coefficients, $a(\zeta)$ and
$b(\zeta)$, are obtained from the asymptotic behavior 
$\vs{\phi}(t;\zeta)\rightarrow (a(\zeta)e^{-i\zeta t}, b(\zeta)e^{i\zeta
t})^{\tp}$ as $t\rightarrow\infty$. The process
of computing these scattering coefficients will be referred 
to as \emph{forward scattering}. 

In general, the nonlinear Fourier spectrum for the potential $q(t)$ comprises 
a \emph{discrete} and a \emph{continuous} spectrum. The discrete spectrum consists 
of the so called \emph{eigenvalues} $\zeta_k\in\field{C}_+$, such that 
$a(\zeta_k)=0$, and, the \emph{norming constants} $b_k$ such that 
$\vs{\phi}(t;\zeta_k)=b_k\vs{\psi}(t;\zeta_k)$. Note that $(\zeta_k,\,b_k)$
describes a \emph{bound state} or a \emph{solitonic state}
associated with the potential. For convenience, let the
discrete spectrum be denoted by the set
\begin{equation}
\mathfrak{S}_K=\{(\zeta_k,\,b_k)\in\field{C}^2|\,\Im{\zeta_k}>0,\,k=1,2,\ldots,K\}.
\end{equation}
The continuous spectrum, also 
referred to as the \emph{reflection coefficient}, is 
defined by $\rho(\xi)={b(\xi)}/{a(\xi)}$ for $\xi\in\field{R}$. In preparation
for the discussion in the following sections, let us define
\begin{equation}\label{eq:a-S}
a_S(\zeta)=\prod_{k=1}^{K}\left(\frac{\zeta-{\zeta}_k}{\zeta-{\zeta}^*_k}\right),
\end{equation}
and $\rho_R(\xi)=a_S(\xi)\rho(\xi)$. The reflection coefficient $\rho_R(\xi)$
now corresponds to a purely radiative potential.

Next, let us note that the class of integrable nonlinear evolution problems that can be treated by the
methods proposed in this article are those described by the Ablowitz-Kaup-Newell-Segur 
formalism~\cite{AKNS1974,AS1981}. In optical fiber communication, the propagation of optical 
field in a loss-less single mode fiber under Kerr-type focusing nonlinearity 
is governed by the nonlinear Schr\"odinger equation (NSE)~\cite{HK1987,Agrawal2013} which 
can be cast into the following standard form
\begin{equation}\label{eq:NSE}
    i\partial_xq=\partial_t^2q+2|q|^2q,\quad(t,x)\in\field{R}\times\field{R}_+,
\end{equation}
where $q(t,x)$ is a complex valued function associated with the slowly varying
envelope of the electric field, $t$ is the retarded time and $x$ 
is position along the fiber. If the potential evolves according
to~\eqref{eq:NSE}, then, the scattering data evolves as:
$b_k(x)=b_ke^{-4i\zeta_k^2x}$ and $\rho(\xi,x)=\rho(\xi)e^{-4i\xi^2x}$ 
($a(\zeta)$ and, consequently, $\zeta_k$ do not evolve). In the rest of 
the paper, we suppress the dependence on $x$ for the sake brevity.

\section{Discrete Inverse Scattering}
\label{sec:discrete-system}
In order to discuss the discretization scheme, we take an equispaced grid defined 
by $t_n= T_1 + nh,\,\,n=0,1,\ldots,N,$ with $t_{N}=T_2$ where $h$ is the grid spacing.
Define $\ell_-,\ell_+\in\field{R}$ such that $h\ell_-= -T_1$, $h\ell_+= T_2$.
Further, let us define $z=e^{i\zeta h}$. For the potential functions supported
in $[T_1, T_2]$, we set $Q_n=2hq(t_n)$, 
$R_{n}=2hr(t_n)$. In the following, we summarize the discrete framework reported
in~\cite{V2017INFT1} which is based on the trapezoidal rule of integration. 
Setting $\Theta_n=1-Q_nR_n$, the recurrence
relation for the Jost solution reads as
$\vv{v}_{n+1}=z^{-1}M_{n+1}(z^2)\vv{v}_n$, which is referred to as the 
\emph{discrete scattering} problem. Here $M_{n+1}(z^2)$ 
is known as the \emph{transfer matrix} which is given by 
\begin{equation}\label{eq:scatter-TR}
M_{n+1}(z^2)=\frac{z^{-1}}{\Theta_{n+1}}
\begin{pmatrix}
1+z^2Q_{n+1}R_n& z^2Q_{n+1}+Q_n\\
R_{n+1}+z^2R_n & R_{n+1}Q_n + z^2
\end{pmatrix}.
\end{equation}
Note that the transfer matrix approach introduced above is analogous 
to that used to solve wave-propagation problems in dielectric 
layered-media~\cite[Chap.~1]{BW1999}. In particular, from the factorization   
\begin{equation*}
\begin{split}
\vv{v}_{n+1}&=\frac{1}{\Theta_{n+1}}
\begin{pmatrix}
1&Q_{n+1}\\
R_{n+1}& 1
\end{pmatrix}
\begin{pmatrix}
z^{-1}&0\\
0& z
\end{pmatrix}
\begin{pmatrix}
1&Q_{n}\\
R_{n}& 1
\end{pmatrix}\vv{v}_n,
\end{split}
\end{equation*}
it can be inferred that the continuous system in~\eqref{eq:zs-prob} is
approximated by two instantaneous scatterers with ``free-space'' 
propagation between them in each of the layers as shown in Fig.~\ref{fig:sm}.
The error analysis of the discrete system 
presented above is carried out in~\cite{V2017INFT1} where it is shown that 
the global order of convergence is $\bigO{h^2}$ for fixed $\zeta$.

In order to express the discrete approximation to the Jost solutions, let us
define the vector-valued polynomial
\begin{equation}\label{eq:poly-vec}
\vv{P}_n(z)=\begin{pmatrix}
            P^{(n)}_{1}(z)\\
            P^{(n)}_{2}(z)
        \end{pmatrix}
         =\sum_{k=0}^n
            \vv{P}^{(n)}_{k}z^k
         =\sum_{k=0}^n
        \begin{pmatrix}
            P^{(n)}_{1,k}\\
            P^{(n)}_{2,k}
        \end{pmatrix}z^k.
\end{equation}
The Jost solution $\vs{\phi}$ can be written in the form $\vs{\phi}_n =
z^{\ell_-}z^{-n}\vv{P}_n(z^2)$ with the initial condition given by 
$\vs{\phi}_0=z^{\ell_-}(1,0)^{\tp}$ that translate into 
$\vv{P}_0=(1,0)^{\tp}$. The recurrence relation for $\vv{P}_n(z^2)$
takes the form
\begin{equation}\label{eq:poly-scatter}
\vv{P}_{n+1}(z^2)= M_{n+1}(z^2)\vv{P}_n(z^2).
\end{equation}
The discrete system discussed above facilitated the development of a 
fast forward scattering algorithm in~\cite{V2017INFT1}. This relied on the fact that the 
transfer matrices have polynomial entries--a form that is amenable to FFT-based
fast polynomial arithmetic~\cite{Henrici1993}.  

In the following sections, we provide details of the fast inverse NFT algorithm
by first developing the methods needed for inversion of the continuous 
spectrum to compute what can be viewed as a purely radiative potential. The 
general version of the inverse NFT is then developed using
the FDT algorithm presented in~\cite{V2017INFT1}.

\subsection{The layer-peeling algorithm}
\label{sec:discrete-TR-summary}
Borrowing the terminology from the theory of layered dielectric
media~\cite[Chap.~1]{BW1999}, let the interval $[t_n,t_{n+1}]$ correspond to the 
$(n+1)$-th \emph{layer} which is
completely characterized by the transfer matrix $M_{n+1}(z^2)$ (see Fig.~\ref{fig:sm}). The
\emph{discrete forward scattering} consists in ``accumulating'' all the layers to 
form $\vv{P}_N(z^2)$. The problem of recovering the discrete samples of the
scattering potential from the discrete scattering coefficients or
$\vv{P}_N(z^2)$ is referred to as the \emph{discrete inverse
scattering} which is facilitated by the so-called 
\emph{layer-peeling} (LP) algorithm. Starting from the recurrence 
relation~\eqref{eq:poly-scatter}, one LP step consists in using 
$\vv{P}_{n+1}(z^2)$ to retrieve the samples of the potential
needed to compute the transfer matrix 
$\wtilde{M}_{n+1}(z^2)=z^{-2}[M_{n+1}(z^2)]^{-1}$ so that the entire step can be
repeated with $\vv{P}_{n}(z^2)$ until all the samples of the potential are 
recovered. The mathematical details of this algorithm can be found
in~\cite{V2017INFT1}. For the sake of reader's convenience, some of the main results 
are summarized below.

Assume $Q_0=0$. Then the recurrence relation~\eqref{eq:poly-scatter} yields
\begin{equation}
\label{eq:TR-cond}
P^{(n+1)}_{1,0}
=\Theta^{-1}_{n+1}\prod_{k=1}^{n}\biggl(\frac{1+Q_kR_k}{1-Q_kR_k}\biggl)
=\Theta^{-1}_{n+1}\prod_{k=1}^{n}\biggl(\frac{2-\Theta_k}{\Theta_k}\biggl),
\end{equation}
and $\vv{P}^{(n+1)}_{n+1}= 0$. The last relationship follows from the assumption $Q_0=0$. For sufficiently 
small $h$, it is reasonable to assume that 
$1+Q_nR_n>0$ so that $P^{(n)}_{1,0}>0$ (it also implies that 
$|Q_n|=|R_n|<1$). The layer-peeling step consists in computing the samples of
the potential, $R_{n+1}$ and $R_n$ (with $Q_{n+1}=-R^*_{n+1}$ and
$Q_{n}=-R^*_{n}$) as follows:
\begin{equation}
R_{n+1} = \frac{P^{(n+1)}_{2,0}}{P^{(n+1)}_{1,0}},\quad
R_n = \frac{\chi}{1 + \sqrt{1+|\chi|^2}},
\end{equation}
where
\[
\chi=\frac{[P^{(n+1)}_{2,1}-R_{n+1}P^{(n+1)}_{1,1}]}{[P^{(n+1)}_{1,0}-Q_{n+1}P^{(n+1)}_{2,0}]}.
\]
Note that $P^{(n+1)}_{1,0}\neq0$ and ${P^{(n+1)}_{1,0} -
Q_{n+1}P^{(n+1)}_{2,0}}\neq0$.
As evident from~\eqref{eq:scatter-TR}, the transfer matrix, $M_{n+1}(z^2)$, 
connecting $\vv{P}_n(z^2)$ and $\vv{P}_{n+1}(z^2)$ is completely determined by
these relations.

If the steps of the LP algorithm are carried out sequentially, one ends up with a
complexity of $\bigO{N^2}$. It turns out that a fast implementation of this LP
algorithm does exist~\cite{V2017INFT1}, which has a complexity of 
$\bigO{N\log^2N}$ for the discrete system 
considered in this article. In the following sections, we describe how to 
synthesize the input for the LP algorithm in order to compute the radiative part of the
scattering potential.



\subsection{Nonlinearly bandlimited signals}
\label{sec:nbs-rho}
A signal is said to be \emph{nonlinearly bandlimited} if it has an empty discrete
spectrum and a reflection coefficient $\rho(\xi)$ that is compactly supported in 
$\field{R}$. This is a direct generalization of the notion of bandlimited signals 
for conventional Fourier transform. However, nonlinearly 
bandlimited signals are not bandlimited, in general. Let us consider 
the reflection coefficient $\rho(\xi)$ as
input. Let the support of $\rho(\xi)$ be contained in $[-\Lambda, \Lambda]$ so
that its Fourier series representation is
\begin{equation}
\rho(\xi)=\sum_{k\in\field{Z}}\rho_k e^{\frac{ik\pi \xi}{\Lambda}}.
\end{equation} 
If $|\rho_k|$ is significant only for $k\geq -n$ ($n\in\field{Z}_+$), then
$\rho(\xi)=\sum_{k=-n}^{\infty}\rho_kz^{2k}+\mathcal{R}_n(z^2)$,
where $z=\exp(i\pi\xi/2\Lambda)$ and $\mathcal{R}_n$ denotes the remainder terms. 
Putting $h=\pi/2\Lambda$ and $T_2=nh\equiv h\ell_+$, we have $\exp(2i\xi T_2)=z^{2n}$ so that
\begin{equation}\label{eq:series}
\breve{\rho}(\xi)=\rho(\xi)z^{2n}
=\sum_{k=0}^{\infty}\breve{\rho}_kz^{2k}+z^{2n}\mathcal{R}_n(z^2).
\end{equation}
Now, it follows that $\breve{\rho}_k=2h\breve{p}(2hk)$ where
\begin{equation}\label{eq:lubich-nbl}
\breve{p}(\tau) = \fourier^{-1}[\rho](\tau)
=\frac{1}{2\pi}\int_{-\Lambda}^{\Lambda}\breve{\rho}(\xi) e^{-i\xi\tau}d\xi.
\end{equation}
Let $2\Lambda_0$ be the fundamental period
and $\Lambda = m\Lambda_0$, where $m\in\field{Z}_+$; then, $h =
\pi/2m\Lambda_0\equiv h_0/m$; therefore, $h\leq h_0$. Now, if we ignore the
remainder term and truncate the series after $N$ terms in~\eqref{eq:series}, 
the input to the fast LP algorithm can be
\begin{equation}
P^{(N)}_1(z^2) = 1,\quad P^{(N)}_2(z^2) =\sum_{k=0}^{N-1}\breve{\rho}_kz^{2k}.
\end{equation} 
This accomplishes the inversion of the reflection coefficient which is assumed
to be compactly supported. Let $\xi_j=j\Delta\xi$ for $j=-M,\ldots,M-1$, where
\[
\Delta\xi = \frac{\pi}{2Mh}.
\]
Then the coefficients $\breve{\rho}_k$ can be estimated using the Fourier sum
\begin{align*}
2hk\breve{p}(2hk) 
&\approx \frac{1}{2M}\sum_{j=-M}^{M-1}\breve{\rho}(\xi_j)e^{-i2hk\xi_j}\\
&= \frac{1}{2M}\sum_{j=-M}^{M-1}\breve{\rho}(\xi_j)e^{-i\frac{2\pi jk}{2M}},
\end{align*}
for $k=0,1,\ldots,N$. The quantity $M$ is chosen to be some multiple of $N$, say,
$M=n_{\text{os}}\times N$ where $n_{\text{os}}$ is referred to as the
\emph{oversampling factor}. Therefore, the overall complexity of synthesizing the
input for the LP algorithm works out to be $\bigO{N\log N}$.

Before we conclude this discussion, let us consider the problem of estimation of
$T_2$. It is of interest to determine a $T_2$ such that the energy in the tail
of the scattering potential, which is to be neglected, is below a certain 
threshold, say, $\epsilon$. Fortunately, there 
is an interesting result due to 
Epstein~\cite{Epstein2004} that allows us to do exactly that. From the theory of Gelfand-Levitan-Marchenko 
equations, it can be shown that there exists a time $T$ such that
\begin{equation}\label{eq:epstien-thm}
\mathcal{E}_+(T)=\int^{\infty}_{T}|q(t)|^2dt
\leq\frac{2\mathcal{I}^2_2(T)}{[1-\mathcal{I}^2_1(T)]},
\end{equation}
assuming $\mathcal{I}_1(T)<1$ where 
\[
\mathcal{I}_m(T)=\left[\int^{\infty}_{2T}|p(-\tau)|^md\tau\right]^{1/m}
\]
for $m=1,2$ (see Appendix~\ref{app:epstein} for a proof which, in essence, is contained 
in the work of Epstein~\cite{Epstein2004}). Let $T=T(\epsilon)$ be such that 
\begin{equation}\label{eq:T-eps}
\frac{2\mathcal{I}^2_2(T)}{[1-\mathcal{I}^2_1(T)]}\leq\epsilon,
\end{equation}
then $\mathcal{E}_+(T)\leq\epsilon$. Consequently, it suffices to choose $T_2\geq
T(\epsilon)$.

\subsubsection{Alternative approach}
\label{sec:nbs-b} 
It is possible to compute the polynomial approximation to the scattering
coefficients $a(\xi)$ and $b(\xi)$ using $\rho(\xi)$, which can be then used to
synthesize the input to the fast LP algorithm. There is no apparent
benefit of this approach compared to the method described above; however, we
describe it for the sake of completeness. The first step consists of constructing 
a polynomial approximation for
$a(\zeta)$ in $|z|<1$ where $z=e^{i\zeta h}$ (under the assumption that no bound
states are present). To this end, let
\begin{equation}\label{eq:b-series}
\rho(\xi)=\sum_{k\in\field{Z}}\rho_kz^{2k},\quad z=e^{i\xi h}.
\end{equation}
With a slight abuse of notation, let us denote this expansion as $\rho(z^2)$.
Let us note that in this case, $a(\xi)$ is not analytic in $\field{R}$ which
means that it is also not analytic on the unit circle $|z|=1$. Here, the
relation~\cite{AKNS1974,AS1981} $|a(\xi)|^2+|b(\xi)|^2=1$ allows us to set up a Riemann-Hilbert (RH) 
problem for a sectionally analytic function 
\begin{equation}
\tilde{g}(z^2)=\begin{cases}
g(z^2) & |z|<1,\\
-{g}^*(1/z^{*2}) & |z|>1,
\end{cases}
\end{equation}
such that the jump condition is given by
\begin{equation}\label{eq:RH-circ}
\tilde{g}^{(-)}(z^2) - \tilde{g}^{(+)}(z^2) =
\log\left[\frac{|\rho(z^2)|^2}{1+|\rho(z^2)|^2}\right],\quad |z|=1,
\end{equation}
where $\tilde{g}^{(-)}(z^2)$ and $\tilde{g}^{(+)}(z^2)$ denotes the boundary values when
approaching the unit circle from $|z|<1$ and $|z|>1$, respectively. Let 
the jump function on the RHS of~\eqref{eq:RH-circ} be denoted by $f(z^2)$ which
can be expanded as a Fourier series
\begin{equation}\label{eq:f-series}
f(z^2)=\sum_{k\in\field{Z}}f_kz^{2k},\quad |z|=1.
\end{equation}
Now, the solution to the RH problem can be stated using the Cauchy integral~\cite[Chap.~14]{Henrici1993} 
\begin{equation}
\tilde{g}(z^2) = \frac{1}{2\pi i}\oint_{|w|=1}\frac{f(w)}{z^2-w}dw.
\end{equation}
The function $g(z^2)$ analytic in $|z|<1$ then works out to be
\begin{equation}
g(z^2)=\sum_{k\in\field{Z}_+\cup\{0\}}f_kz^{2k},\quad |z|<1.
\end{equation}
Finally, $a_N(z^2)=\{\exp[g(z^2)]\}_{N}$ with $z=e^{i\zeta h}$ where
$\{\cdot\}_N$ denotes truncation after $N$ terms. The implementation of the 
procedure laid out above can be carried out using the FFT
algorithm, which involves computation of the coefficients $f_k$ and the 
exponentiation in the last step~\cite[Chap.~13]{Henrici1993}. Note that, in the computation of
$g(z^2)$, we discarded half of the coefficients; therefore, in the numerical implementation
it is necessary to work with at least $2N$ number of samples of
$f(z^2)$ in order to obtain $a_N(z^2)$ which is a polynomial of degree $N-1$.

The next step is to compute the polynomial approximation for $\breve{b}(\xi)$. To
this end, consider
\begin{equation}
\breve{b}(\xi)=b(\xi)z^{2n}
=\left[\sum_{k=0}^{\infty}\breve{\rho}_kz^{2k}+z^{2n}\mathcal{R}_n(z^2)\right]\exp[g(z^2)].
\end{equation}
In the following, we will again discard the remainder term. The polynomial
approximation for $\breve{b}(\xi)$ reads as
\begin{equation}
\breve{b}_N(z^2)=\left\{a_N(z^2)\sum_{k=0}^{N-1}\breve{\rho}_kz^{2k}\right\}_N
=\sum_{k=0}^{N-1}\breve{b}_kz^{2k},
\end{equation}
Now, the input to the fast LP algorithm works out to be
\begin{equation}
P^{(N)}_1(z^2)=\sum_{k=0}^{N-1}a_kz^{2k},\quad P^{(N)}_2(z^2)=\sum_{k=0}^{N-1}\breve{b}_kz^{2k}.
\end{equation}

\subsection{Fast inverse NFT}
\label{sec:fast-inverse-NFT}

In the previous sections, we restricted ourselves to the case of empty discrete
spectrum. In this section, we describe how a fast inverse NFT algorithm can be 
developed for the general NF spectrum 
using either the Classical DT (CDT) or the FDT algorithm reported 
in~\cite{V2017INFT1}. 
Given a reflection coefficient $\rho(\xi),\,\xi\in\field{R},$ and
the discrete spectrum $\mathfrak{S}_K$, define 
$a_S(\xi)$ as in~\eqref{eq:a-S} and $\rho_R(\xi)=a_S(\xi)\rho(\xi)$. Let $q(t)$ 
denote the scattering potential
corresponding to the aforementioned NF spectrum. 

Now, as illustrated in Fig.~\ref{fig:schema-inft}, the inverse NFT can be carried out 
in the following two steps:
\begin{itemize}
\item[I.] Generate the signal $q_R(t)$ corresponding to the reflection
coefficient $\rho_R(\xi)$ using the method described in 
Sec.~\ref{sec:nbs-rho}. This amounts to computing the purely radiative part of 
the complete potential $q(t)$. The complexity of this step is 
$\bigO{N\log^2N}$ if the number of nodes used for the FFT operation
involved there is given by $M=n_{\text{os}}N$ where $n_{\text{os}}\ll N$. Here,
$n_{\text{os}}$ can be identified as the 
oversampling factor (typically $\leq8$).

\item[II.] Use the signal $q_R(t)$ as the seed potential and add
bound states described by $\mathfrak{S}_K$ using the CDT or the FDT algorithm to
obtain $q(t)$. The complexity of this step is $\bigO{N(K+\log^2N)}$
when FDT is employed while $\bigO{K^2N}$ when CDT is employed. Here we
also consider the partial-fraction (PF) variant of the FDT
algorithm (labeled as FDT-PF), which is shown to offer a small increase in speed~\cite{V2017INFT1}. 

Finally, let us note that the overall complexity of the inverse NFT is given by
$\bigO{N(K+\log^2N)}$ when FDT is used and $\bigO{N(K^2+\log^2N)}$ when CDT is
used.
\end{itemize}

\section{Numerical Experiments}
\label{sec:num-res}
Let $q^{(\text{num}.)}$ denote the numerically computed potential for a given NF
spectrum. If the exact potential $q$ is known, then we quantify the error as 
\begin{equation}\label{eq:e_rel-q}
e_{\text{rel.}}={\|q^{(\text{num.})}-q\|_{\fs{L}^2}}/{\|q\|_{\fs{L}^2}},
\end{equation}
where the integrals are evaluated numerically using the trapezoidal rule. For 
the purpose of convergence analysis, only those examples are deemed to be admissible 
where closed-form solutions are available. However,
on account of scarcity of such examples, an exhaustive test for universality of 
the algorithm cannot be carried out in this manner. To remedy this, we choose a
higher-order convergent algorithm for the forward scattering problem and compute the NF
spectrum of the potential generated by the fast inverse NFT. The 
error between the computed NF spectrum and the provided NF spectrum can serve 
as a good metric to measure the robustness of the algorithm. 

For the higher-order 
scheme, we choose the (exponential) $3$-step \emph{implicit Adams} 
method (IA$_3$)~\cite{HNW1993} which has an order of convergence $4$, i.e.,
$\bigO{N^{-4}}$ (see Appendix~\ref{app:IA} for details). Fortunately, 
this method can also be made fast by the use of
FFT-based polynomial arithmetic which allows us to test for large number 
of samples ($N\in\{2^{10},2^{11},\ldots,2^{20}\}$). Note that this procedure 
by no means qualifies as the test for total numerical error 
on account of the fact that the error metric is not the \emph{true} numerical 
error. Therefore, the results in this case must be interpreted with caution.
Further, for the sake of comparison, we also consider the \emph{T\"oplitz inner bordering} 
(TIB) algorithm for inverse scattering 
(Belai~\et~\cite{BFPS2007}) whenever the
discrete spectrum is empty. We use the second order convergent version of this
algorithm which has also been reported in Frumin~\et~\cite{FBPS2015}. The latter
paper provides an improved understanding of the
original algorithm presented in~\cite{BFPS2007}; therefore, we choose to refer 
to~\cite{FBPS2015} in this article whenever we mention the TIB algorithm.

Finally, let us emphasize that the primary objective behind the numerical tests in 
this section is to verify the trends expected from the theory. The actual values
of any defined performance metric observed in the results are merely 
representative of what can be achieved\footnote{The total run-time, for
instance, may differ on different computing machines; therefore, we would only be 
interested in trends as far as the complexity analysis of the 
algorithms are concerned.}, and, admittedly, better results 
can be obtain by appropriately tuning the parameters used in the
test. For instance, a good choice of the computational domain
helps to maintain a smaller step-size in the numerical discretization and,
hence, lowers the numerical error.

\begin{figure*}[!tbh]
\centering
\includegraphics[scale=1]{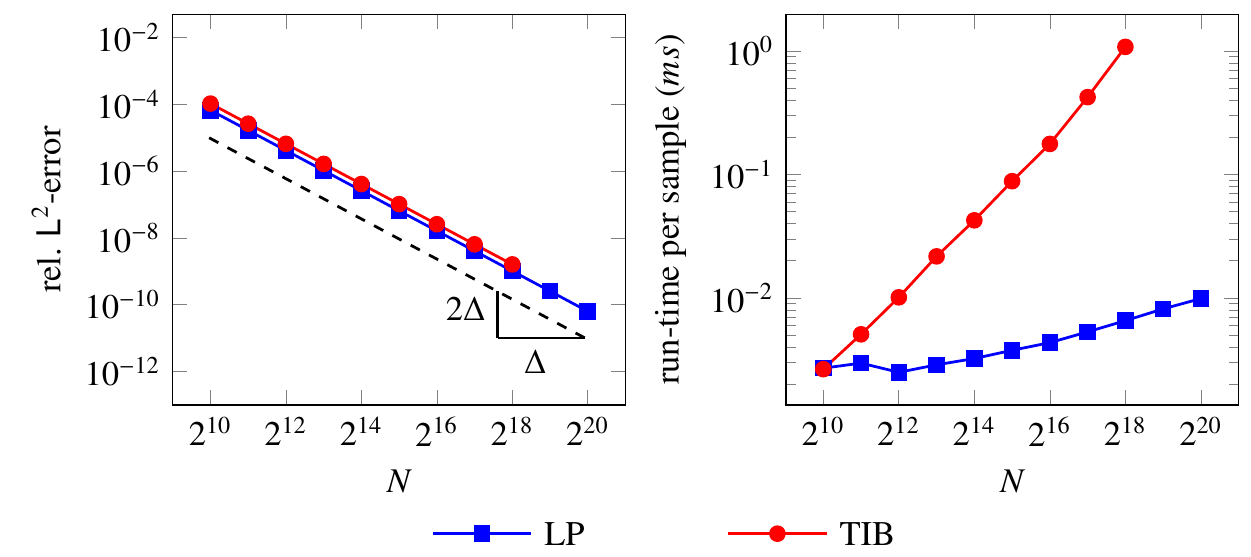}
\caption{\label{fig:TR-comp}The figure shows a comparison of the 
algorithms LP and TIB for the secant-hyperbolic
potential ($A_R=0.4$) with respect to convergence rate (left) and run-time per sample (right).}
\end{figure*}

\begin{figure*}[!th]
\centering
\def\scale{1}
\includegraphics[scale=\scale]{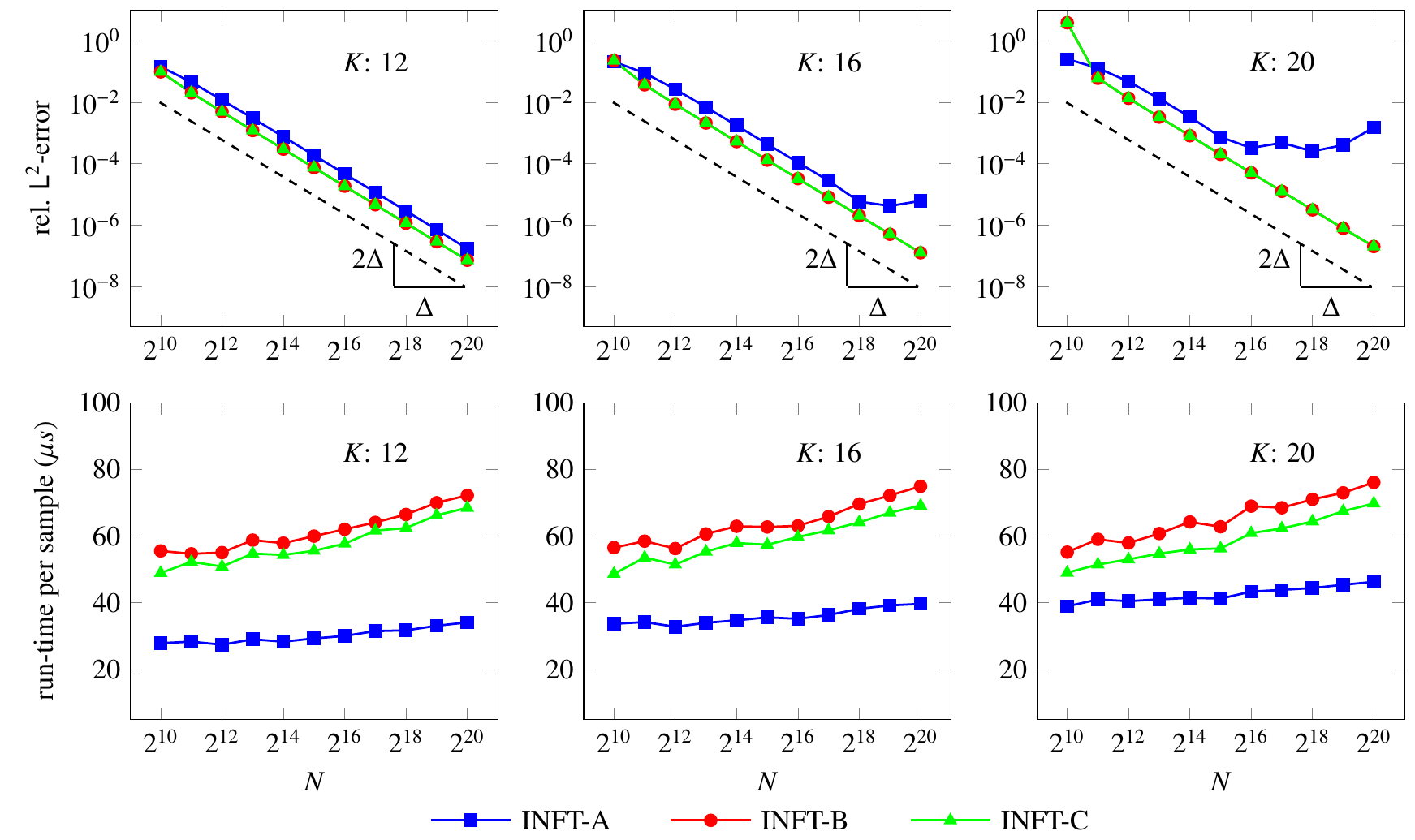}
\caption{\label{fig:fast-INFT}The figure shows the performance of the algorithms
INFT-A/-B/-C for a fixed number of eigenvalues ($K\in\{12,16,20\}$) and
varying number of samples ($N$) for the secant-hyperbolic potential 
(see Sec.~\ref{sec:sech-ex}). The error plotted on the
vertical axis is defined by~\eqref{eq:ds-sech}.}
\end{figure*}

\begin{figure*}[!th]
\centering
\def\scale{1}
\includegraphics[scale=\scale]{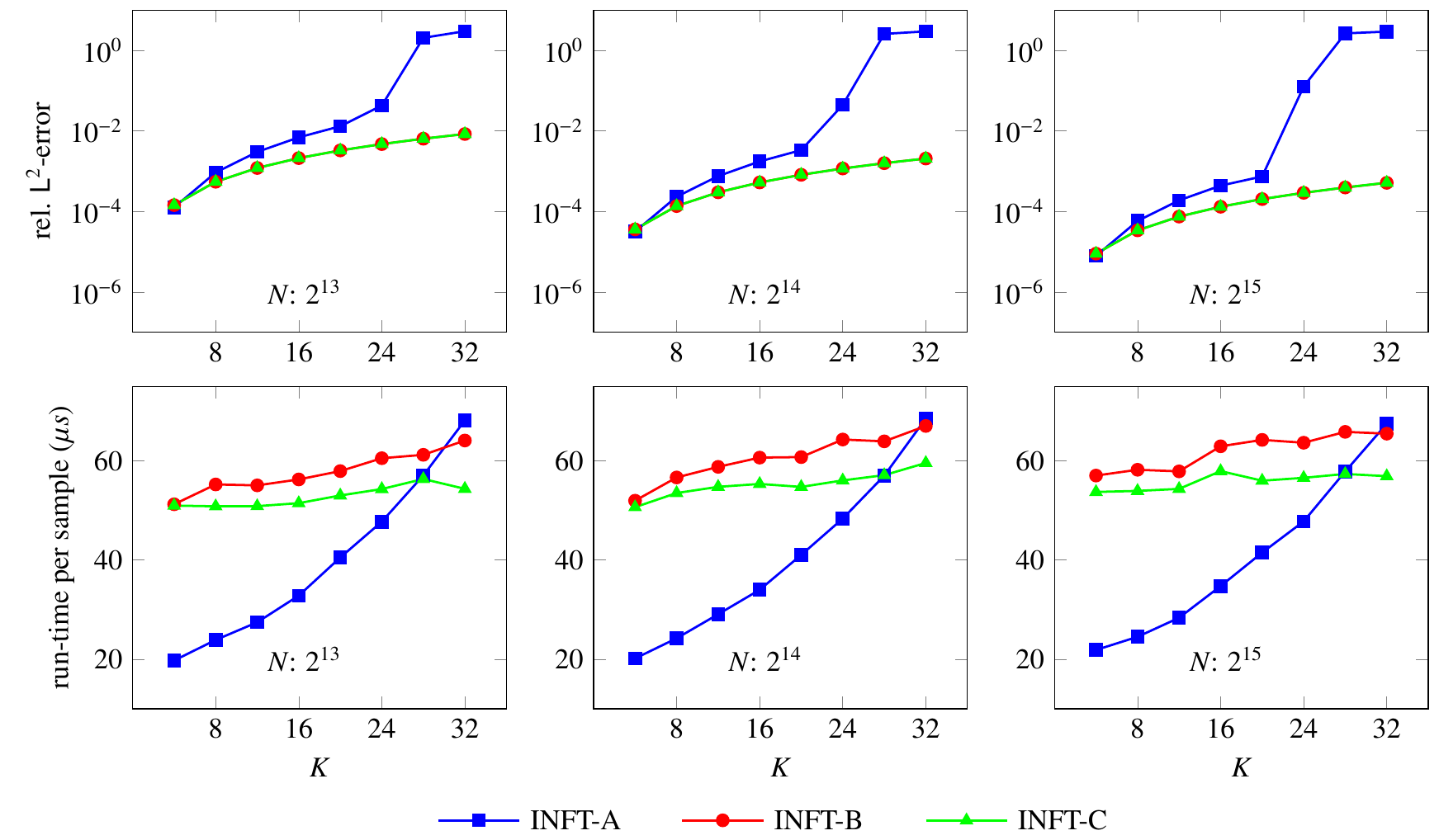}
\caption{\label{fig:fast-INFT-EV}The figure shows the performance of the algorithms 
INFT-A/-B/-C for a fixed number of samples ($N\in\{2^{13},2^{14},2^{15}\}$) and varying number of eigenvalues
($K$) for the secant-hyperbolic potential (see Sec.~\ref{sec:sech-ex}).
The error is quantified by~\eqref{eq:ds-sech}.}
\end{figure*}

\subsection{Secant-hyperbolic potential}
\label{sec:sech-ex}
Here, we would like to devise tests to confirm the order of
convergence and the complexity of computations for the algorithms proposed 
thus far. To this end, we choose the secant-hyperbolic potential given by 
$q(t)=(A_R+K)\sech t$, which is treated exactly in~\cite{SY1974}. Here 
$A_R\in[0,0.5)$ and $K$ is a positive integer. The discrete spectrum can be stated as
\begin{equation}\label{eq:ds-sech}
\mathfrak{S}_{K}=\left\{(\zeta_k,b_k)
\left|\begin{aligned}
&\zeta_k=i(A_R+0.5+K-k),\\ 
&b_k=(-1)^{k},\,k=1,2,\ldots,K\end{aligned}\right.
\right\},
\end{equation}
and the continuous spectrum is given by $\rho=\rho_R/a_S$
where $a_S(\xi)$ is defined by~\eqref{eq:a-S} and
\begin{equation}\label{eq:cs-R}
\rho_R(\xi)=b(\xi)\frac{\Gamma(0.5+A_R-i\xi)\Gamma(0.5-A_R-i\xi)}{[\Gamma(0.5-i\xi)]^2}.
\end{equation}
with $b(\xi)=-\sin[(A_R+K)\pi]\sech(\pi\xi)$. This test consist 
in studying the behavior of the fast INFTs 
for different number of samples ($N$) as well as eigenvalues ($K$).
We set $A_R=0.4$. The scattering potential is scaled by 
$\kappa= 2(\sum_{k=1}^K\Im\zeta_k)^{1/2}$ and $[-T,\,T]$, 
$T=30\kappa/\min_{k}(\Im\zeta_k)$, is taken as the computational
domain and we set 
$N_{\text{th}}=N/8$ for FDT-PF as in~\cite{V2017INFT1}.

Let us first consider the case $K=0$ so that $\rho=\rho_R$ (setting the convention that
$a_S=1$ when $K=0$). Note that on account of the exponential decay of $\rho$, it can 
be assumed to be effectively supported in a bounded domain. Besides the knowledge of the true
potential allows us to provide a good estimate of the computational domain. 
Set $T=\log(2A_R/\epsilon)\approx 30$ for $\epsilon=10^{-12}$, then $[-T,\,T]$ 
can be taken as the computational domain\footnote{For the ZS problem, 
let us note that the error in the initial condition at the left-boundary 
can be kept below $\epsilon>0$, if 
$\|q\chi_{(-\infty,T_1]}\|_{\fs{L}^1}\leq\sinh^{-1}\epsilon$~\cite{V2017INFT1}.}. The 
result for $A_R=0.4$ is plotted in Fig.~\ref{fig:TR-comp} which shows that 
the performance of LP is comparable to that of TIB. Further, each of these algorithms 
exhibit a second order of convergence (i.e., error vanishing as
$\bigO{N^{-2}}$). The run-time 
behavior in Fig.~\ref{fig:TR-comp} shows that LP-based INFTs have a poly-log complexity 
per sample as opposed to the $\bigO{N}$ complexity per sample exhibited by TIB.

For $K>0$, the results are plotted in Fig.~\ref{fig:fast-INFT} which reveal that the fast INFTs 
based on FDT (labeled as `INFT-B') and FDT-PF (labeled as `INFT-C') are superior to 
that based on CDT (labeled as `INFT-A') which becomes unstable with increasing 
number of eigenvalues. The latter, however, can be useful for a small number 
eigenvalues. The figure also confirms the second
order of convergence of INFT-B/-C which is consistent with the underlying 
one-step method, namely, the trapezoidal rule. For small number of eigenvalues, 
INFT-A also exhibits a second order of convergence. Finally, let us 
observe that, for fixed $N$, INFT-A has a 
complexity of $\bigO{K^2}$ and that for INFT-B/-C is $\bigO{K}$. While these trends can be 
confirmed from Fig.~\ref{fig:fast-INFT-EV}, let us mention that, with an improved
implementation, INFT-B/-C can be made even more competitive to INFT-A in complexity. 

\begin{figure}[!th]
\centering
\def\scale{1}
\includegraphics[scale=\scale]{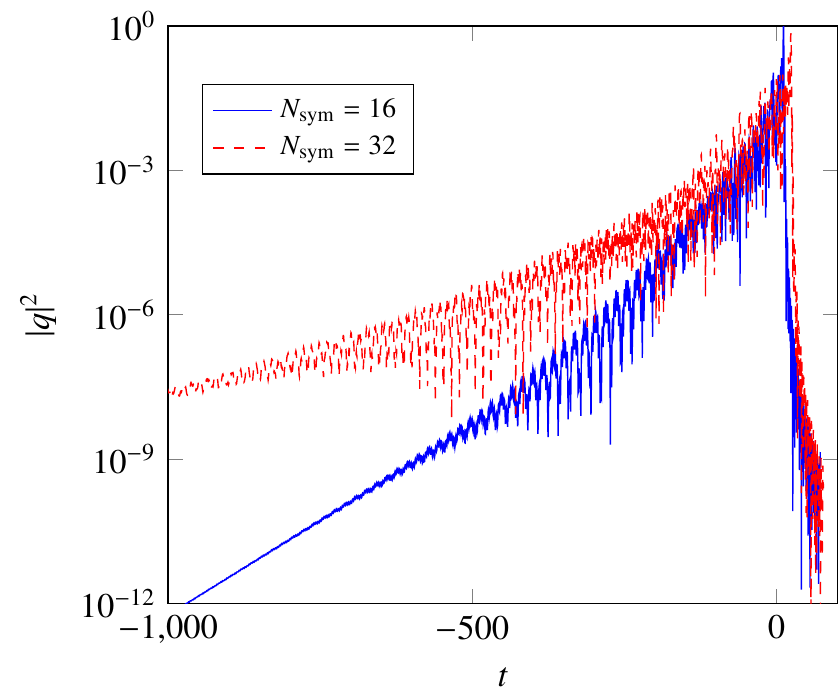}
\caption{\label{fig:QPSK-sig}The figure shows the potential corresponding to a QPSK modulated
continuous spectrum given by~\eqref{eq:QPSK-rho} with number of symbols
$N_{\text{sym}}\in\{16,32\}$. The number of samples used is $N=2^{12}$ and 
the computational domain is $[-15T_2, T_2]$ where $T_2$ is given 
by~\eqref{eq:QPSK-T2}. Also, we set $A_{\text{eff.}}=10$ which is defined
by~\eqref{eq:A_eff}.}
\end{figure}

\begin{figure*}[!th]
\centering
\def\scale{1}
\includegraphics[scale=\scale]{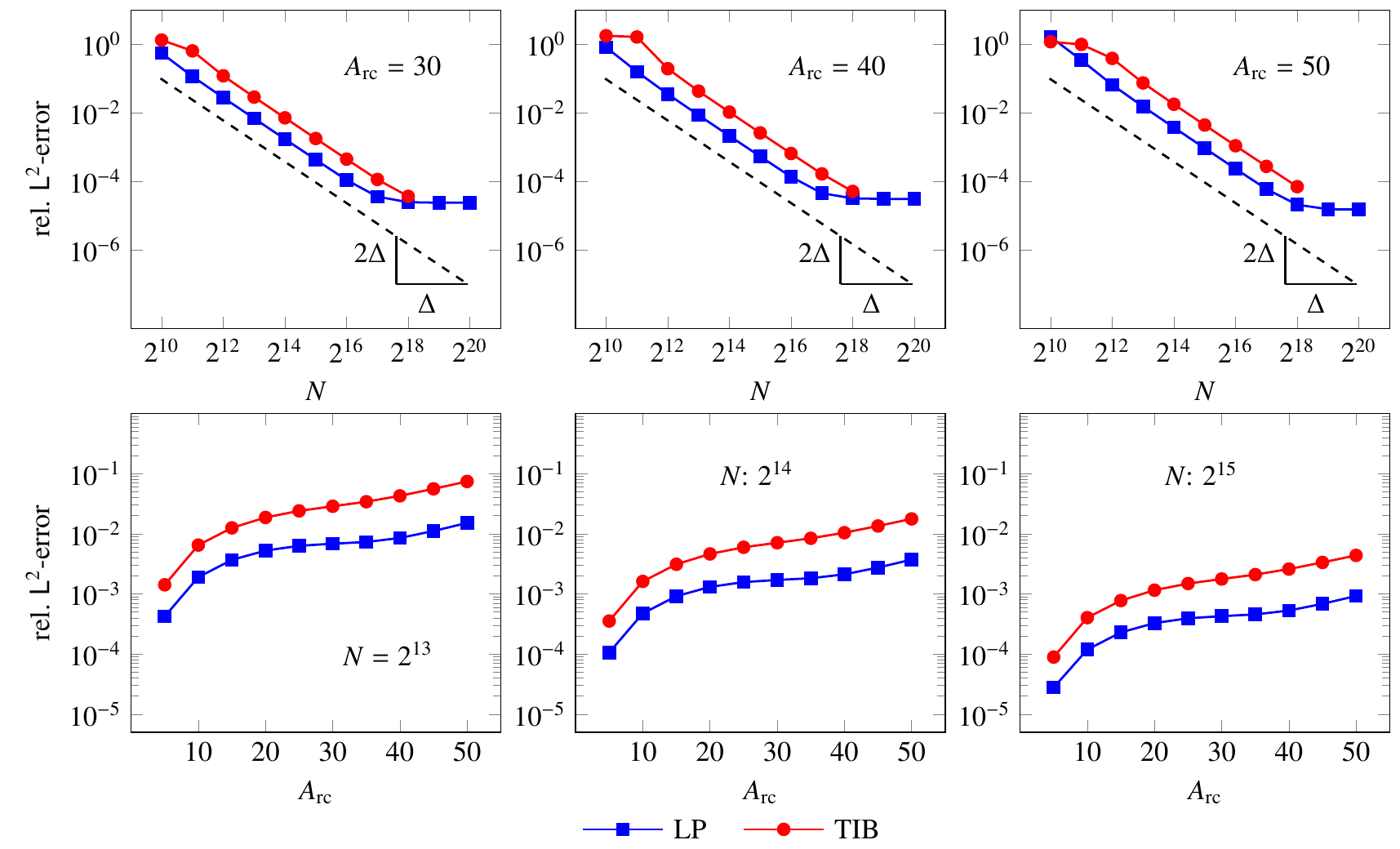}
\caption{\label{fig:convg-rc} The figure shows the error
analysis for the signal generated from the continuous spectrum given 
by~\eqref{eq:RC-spec} which is the frequency-domain description of the raised-cosine
filter (see Sec.~\ref{sec:nbs-test}). The error is quantified
by~\eqref{eq:e_rel-rho}.}
\end{figure*}

\begin{figure*}[!th]
\centering
\def\scale{1}
\includegraphics[scale=\scale]{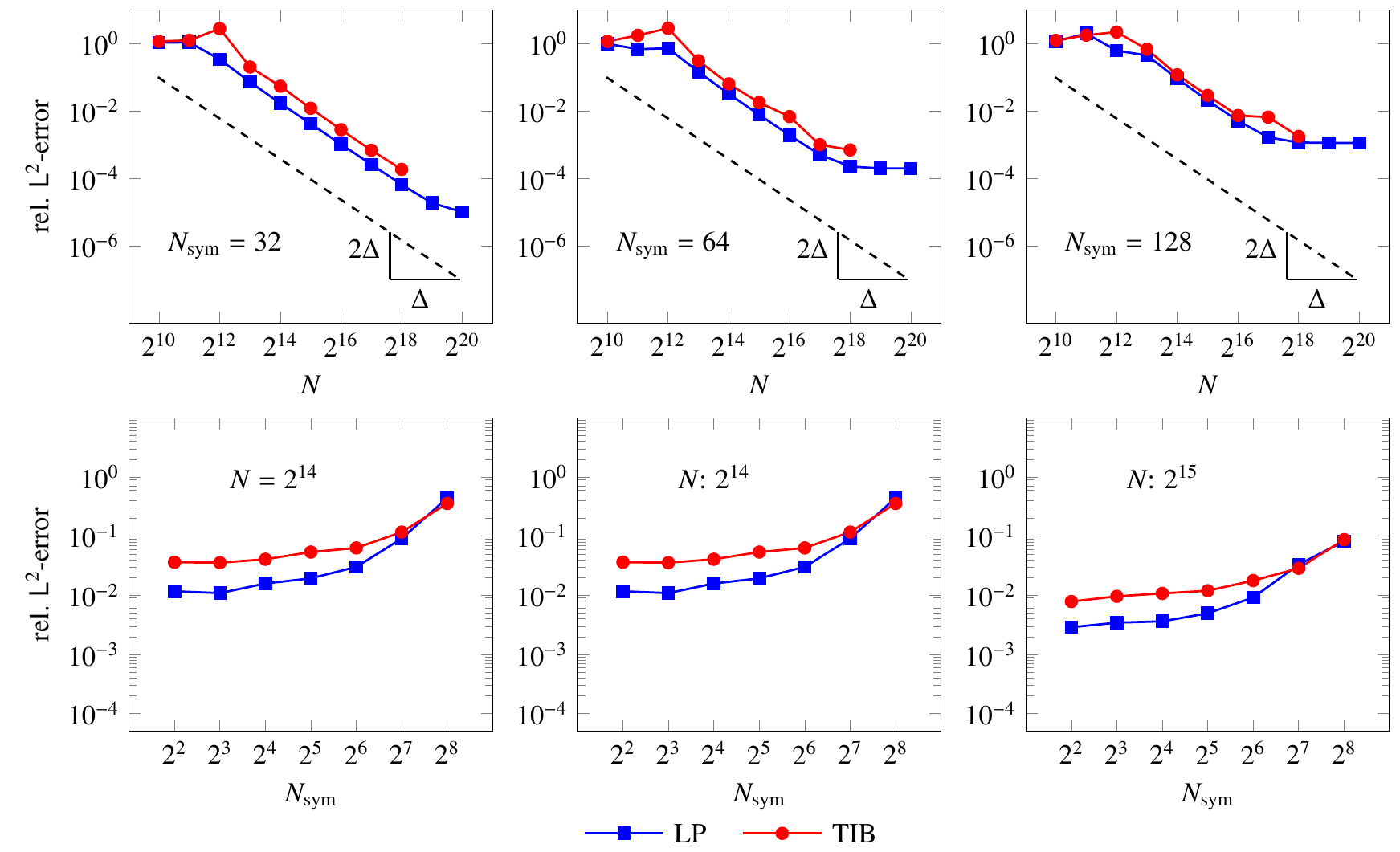}
\caption{\label{fig:convg-qpskrc} The figure shows the error analysis for
the QPSK modulated continuous spectrum given by~\eqref{eq:QPSK-rho} for varying
number of symbols $N_{\text{sym}}$ (see Sec.~\ref{sec:nbs-test}). Here, we 
set $A_{\text{eff.}}=10$ which is defined
by~\eqref{eq:A_eff}.}
\end{figure*}

\subsection{Nonlinearly bandlimited signals}
\label{sec:nbs-test}
Let us consider a soliton-free signal whose continuous spectrum is given by 
\begin{equation}\label{eq:RC-spec}
H_{\text{rc}}(\xi)=
\begin{cases}
A_{\text{rc}}&|\tau_s\xi|\leq 1-\beta,\\
\frac{A_{\text{rc}}}{2}
\left[1+\cos\left(\frac{\pi}{2\beta}\Xi\right)\right]&
||\tau_s\xi|-1|\leq{\beta},\\
0 & |\tau_s\xi|>{1+\beta},
\end{cases}
\end{equation}
where $\Xi=|\tau_s\xi|-(1-\beta)$ with $\beta\in[0,1]$, and, 
$A_{\text{rc}}$ and $\tau_s$ are positive constants. 
The nonlinear impulse response (NIR) $h_{\text{rc}}(\tau)$ can be worked out exactly; 
however, we do not use this information for constructing the input to the fast
LP algorithm. Note that $H_{\text{rc}}(\xi)$ and $h_{\text{rc}}(\tau)$ describe 
the well-known \emph{raised-cosine} filter in the frequency-domain and the time-domain, respectively.

In order to estimate the computational domain, we use Epstein's result discussed
in Sec.~\ref{sec:nbs-rho} which consists in finding a time $T$ such that
\begin{equation}\label{eq:epstien}
\mathcal{E}_+(T)=\int^{\infty}_{T}|q(t)|^2dt\leq\frac{2\mathcal{I}^2_2(T)}{[1-\mathcal{I}^2_1(T)]},
\end{equation}
assuming $\mathcal{I}_1(T)<1$ where 
\[
\mathcal{I}_m(T)=\left[\int^{\infty}_{2T}|h_{\text{rc}}(-\tau)|^md\tau\right]^{1/m}
\]
for $m=1,2$. A crude estimate for $T$ such that 
$\mathcal{I}^2_2(T)=\epsilon$ is given by
\begin{equation}
2T(\epsilon)\sim\left({A_{\text{rc}}^2\tau_s^4}\right)^{1/5}
\beta^{-4/5}\epsilon^{-1/5},
\end{equation}
which uses the asymptotic form of $h_{\text{rc}}(\tau)$. If $\epsilon\ll1$, we
may assume that the potential is effectively supported\footnote{The Epstein's theorem 
provides an estimate for the right boundary if the right NIR
is used; therefore, strictly speaking, the computational domain must be of the 
form $(-\infty, T(\epsilon)]$.} in
$[-T(\epsilon),T(\epsilon)]$ where we set $\epsilon=10^{-9}$. Also, let
$\beta=0.5$ and $\tau_s=1$ in the following. 

For this example, we devise two kinds of tests. For the first kind of tests, we
disregard any modulation scheme and carry out the inverse NFT for varying number
of samples ($N$) for each of the values of $A_{\text{rc}}\in\{10,\ldots,50\}$.
In the second kind of tests, we considers the quadrature-phase-shift-keyed
(QPSK) modulation scheme which is
described later. Let $\Omega_h = [-\pi/2h,\pi/2h]$, then 
the error is quantified by 
\begin{equation}\label{eq:e_rel-rho}
e_{\text{rel.}}=
{\|\rho^{(\text{num.})}-\rho\|_{\fs{L}^2(\Omega_h)}}/{\|\rho\|_{\fs{L}^2(\Omega_h)}},
\end{equation}
where the integrals are computed from $N$ equispaced samples in $\Omega_h$ using the
trapezoidal rule. As stated in the beginning, the quantity $\rho^{(\text{num.})}$
is computed using the (exponential) IA$_3$.

The results of the first kind of tests are shown in Fig.~\ref{fig:convg-rc} where a
comparison is made between LP and 
TIB\footnote{The complexity of TIB becomes prohibitive for increasing
$N$, therefore, we restrict ourselves to $N\leq2^{18}$.}. From the plots in the
top row of Fig.~\ref{fig:convg-rc}, the second order of
convergence is readily confirmed for both of these algorithms with LP performing
somewhat better than TIB. The plateauing of the error in these plots can be
attributed to accumulating numerical errors in the inverse NFT algorithm as well as
the implicit Adams method. The behavior of the error with respect to 
$A_{\text{rc}}$ is shown in the bottom row of Fig.~\ref{fig:convg-rc} where
LP shows better performance than TIB. 

Now, for the second kind of tests, we consider the QPSK 
modulation of the continuous spectrum as follows
\begin{equation}\label{eq:QPSK-rho}
\rho(\xi)=\left(\sum_{n\in J}
s_ne^{-in\pi\tau_s\xi}\right)H_{\text{rc}}(\xi)=S(\xi)H_{\text{rc}}(\xi),
\end{equation}
where the index set is 
$J=\{-N_{\text{sym}}/2,\ldots,N_{\text{sym}}/2-1\}$ and $s_n\in\{\pm1, \pm i\}$
with $N_{\text{sym}}>0$ being an even integer. The 
estimate for the right boundary works out to be 
\begin{equation}\label{eq:QPSK-T2}
T_2=T(\epsilon)+\pi\tau_s N_{\text{sym}}/4;
\end{equation}
however, an estimate for the for the left boundary is not available in a closed
form. Here, we take a heuristic approach by setting $T_1=-W \times T_2$ where $W$
is chosen by trial and error. The scale factor $A_{\text{rc}}$ is chosen
such that $A_{\text{eff.}}=10$ where
\begin{equation}\label{eq:A_eff}
A_{\text{eff.}}=\|\rho\|_2/\|H_{\text{rc}}\|_2. 
\end{equation}
It is important
to observe here that the signal generated from~\eqref{eq:QPSK-rho} is highly 
asymmetric with poor decay behavior as
$t\rightarrow-\infty$ (see Fig.~\ref{fig:QPSK-sig}). The higher values of the quantities 
$N_{\text{sym}}$ and $A_{\text{eff.}}$, both, tend to worsen this phenomenon.
Therefore, this example turns out to be very challenging for the numerical
algorithm. In Fig.~\ref{fig:convg-qpskrc}, we provide results of numerical experiments conducted with 
$N_{\text{sym}}\in\{4,8,\ldots,256\}$ number of symbols where
$W=5\log_2N_{\text{sym}}$ is used to determine the computational domain. The accuracy of 
LP and TIB, both, tends to worsen with increasing number of symbols where 
LP performs slightly better than TIB. Based on these results it is evident
that any method of pulse-shaping must take into account the
relationship between the signal and its NF spectrum as opposed to directly
applying conventional Fourier transform based techniques of
pulse-shaping.

\begin{figure*}[!th]
\centering
\def\scale{1}
\includegraphics[scale=\scale]{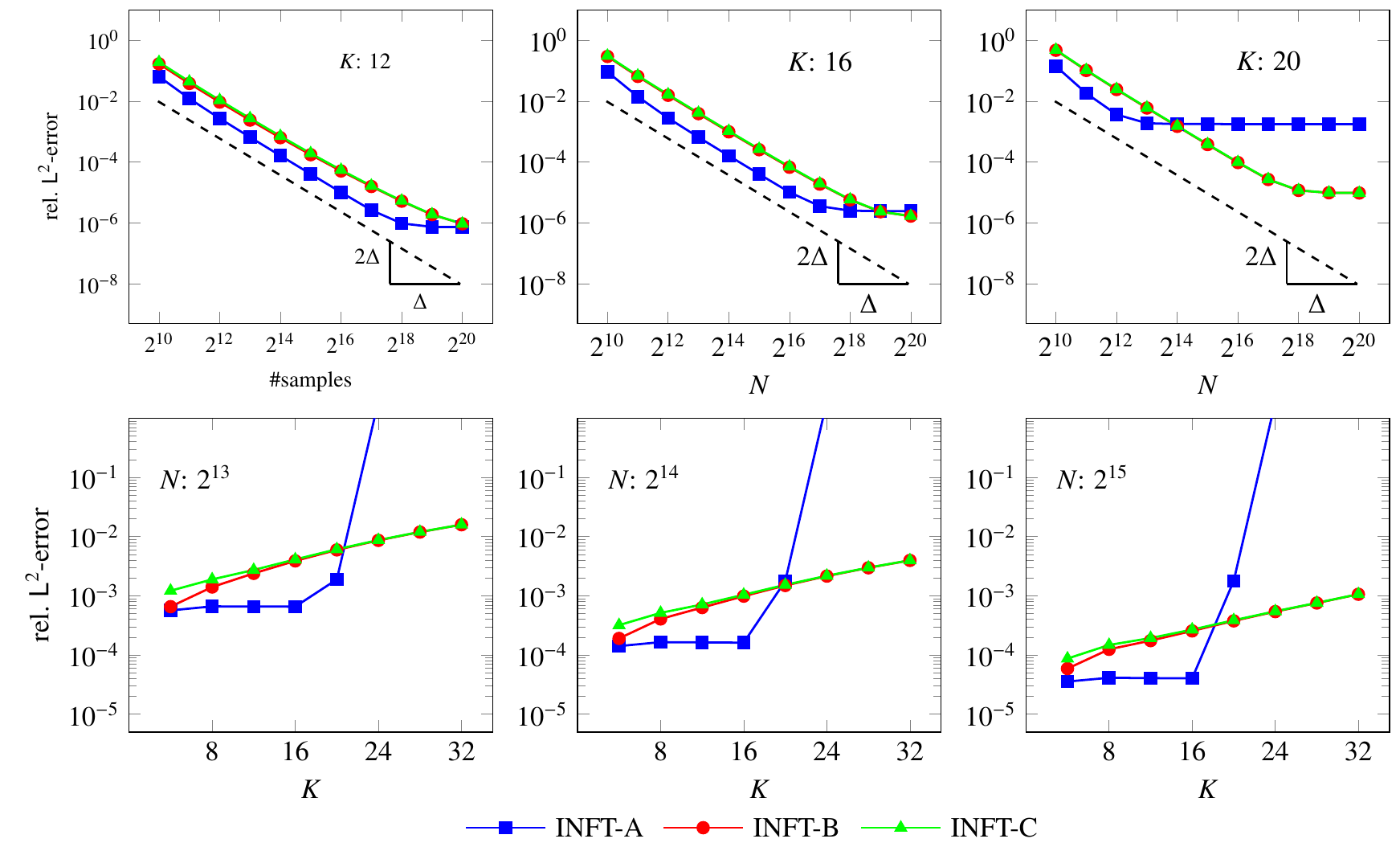}
\caption{\label{fig:convg-rc-dt} The figure shows the results of the error
analysis for an example where the discrete spectrum is $\mathfrak{S}_K$ and the continuous spectrum is
identical to the Fourier spectrum of the raised-cosine filter (see
Sec.~\ref{sec:bound-states-rc}). The error plotted on the vertical axis corresponds to
the continuous spectrum, which is quantified by~\eqref{eq:e_rel-rho}. Here
$A_{\text{rc}}=20$.}
\end{figure*}

\subsubsection{Addition of Bound states}
\label{sec:bound-states-rc}
Here, we fix $A_{\text{rc}}=20$ and assume no modulation of the continuous spectrum. The bound states 
to be added are described by~\eqref{eq:ds-sech}.  Let us observe that the ``augmented''
potential has a reflection coefficient which is given by
$\rho^{(\text{aug.})} = \rho/a_S$. Now the error can be 
quantified by~\eqref{eq:e_rel-rho}. The potentials are scaled by $\kappa$ as in
Sec.~\ref{sec:sech-ex} and the computational domain is chosen such that
$-T_1=T_2=T(\epsilon)\kappa/\min_{k}(\Im\zeta_k)$.

The results for the continuous spectrum are shown in Fig.~\ref{fig:convg-rc-dt}
where the order of convergence can be confirmed from the plots in the top
row. The plots in the bottom row reveal that INFT-A, which is based on CDT, is 
unstable for increasing number of eigenvalues. On the other hand,
the algorithms INFT-B/-C, which are based on FDT/FDT-PF, respectively, seem to 
perform equally well without showing any signs of instability.

\begin{figure*}[!th]
\centering
\def\scale{1}
\includegraphics[scale=\scale]{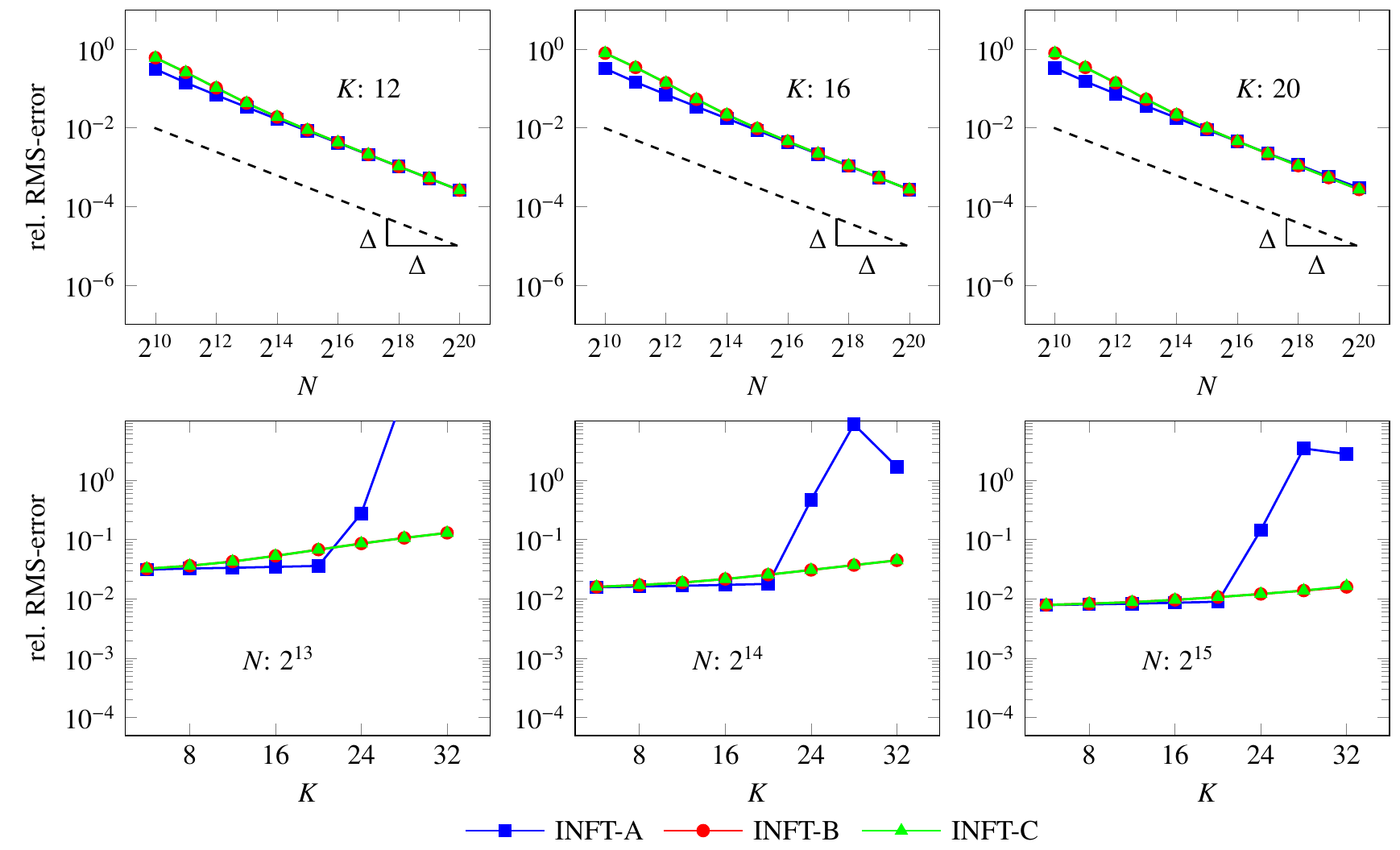}
\caption{\label{fig:convg-rc-dt-nconst} The figures shows the results of the error
analysis for an example where the discrete spectrum is $\mathfrak{S}_K$ and the continuous spectrum is
identical to the Fourier spectrum of the raised-cosine filter (see
Sec.~\ref{sec:nbs-test}). The error plotted on the vertical axis corresponds to
the norming constant, which is quantified by~\eqref{eq:err-RMS}. Here
$A_{\text{rc}}=20$.}
\end{figure*}

For the discrete spectrum, we assume that the discrete eigenvalues are known
exactly, and, then use this information to compute the norming constants using
the method discussed in~\cite{V2017INFT1}. The
error is quantified by
\begin{equation}\label{eq:err-RMS}
e_{\text{rel.}}=\sqrt{\left({\sum_{k=1}^K|b^{(\text{num.})}_k-b_k|^2}\right)\biggl/
{\sum_{k=1}^K|b_k|^2}},
\end{equation}
where $b^{(\text{num.})}_k$ is the numerically computed norming constant using IA$_3$.
The results are shown in Fig.~\ref{fig:convg-rc-dt-nconst}
where the order of convergence turns out to be $\bigO{N^{-1}}$ from the plots 
in the top row. This decrease of
order of convergence can be attributed to the use of the true eigenvalues 
as opposed to the numerically computed one to compute the norming constants. Again, 
the plots in the bottom row reveal that INFT-A is unstable for increasing number 
of eigenvalues. On the other hand, the algorithms INFT-B/-C seem to perform equally well while 
showing no signs of instability.

\section{Conclusion}
\label{sec:final}
To conclude, we have presented two new fast INFT algorithms 
with $\bigO{KN+N\log^2N}$ complexity and a 
convergence rate of $\bigO{N^{-2}}$. These algorithms are based on the discrete 
framework introduced in~\cite{V2017INFT1} for the ZS scattering problem 
where the well-known one-step method, namely, the 
\emph{trapezoidal rule} is employed for the numerical discretization. Further, our 
algorithm depends on the fast LP and the FDT algorithm presented 
in~\cite{V2017INFT1}. Numerical tests reveal 
that both variants of the INFT algorithm are 
capable of dealing with a large number of eigenvalues (within the limitations of
the double precision arithmetic) previously unreported.
Further, for the cases considered in this article, our algorithms perform
better than the TIB algorithm~\cite{BFPS2007,FBPS2015} in terms of 
accuracy while being faster by an order of magnitude. Let us also note that the
TIB algorithm has no consequence for the fast inverse NFT in the general case.

Next, let us mention that 
we have not included simulations of a realistic 
optical fiber link in order to demonstrate the effectiveness of our algorithms. A thorough 
testing for various NFT-based modulation schemes for a realistic optical fiber 
link is beyond the scope of this paper. This omission however does not impact
the study of the limitation of the proposed algorithms from a numerical analysis
perspective. 

Future research on fast INFTs will further focus on the stability properties of the
LP algorithm and the DT iterations. Moreover, 
we would also like to consider other linear multistep
methods to obtain a higher-order convergent forward/inverse NFTs. The implicit
Adams method used in this paper for the purpose of testing already demonstrates
that such possibilities do exist, at least, for the solution of the direct ZS
problem.



\providecommand{\noopsort}[1]{}\providecommand{\singleletter}[1]{#1}%
%


\appendix

\section{Implicit Adams Method}
\label{app:IA}

\begin{figure*}[ht!]
\centering
\includegraphics[scale=1]{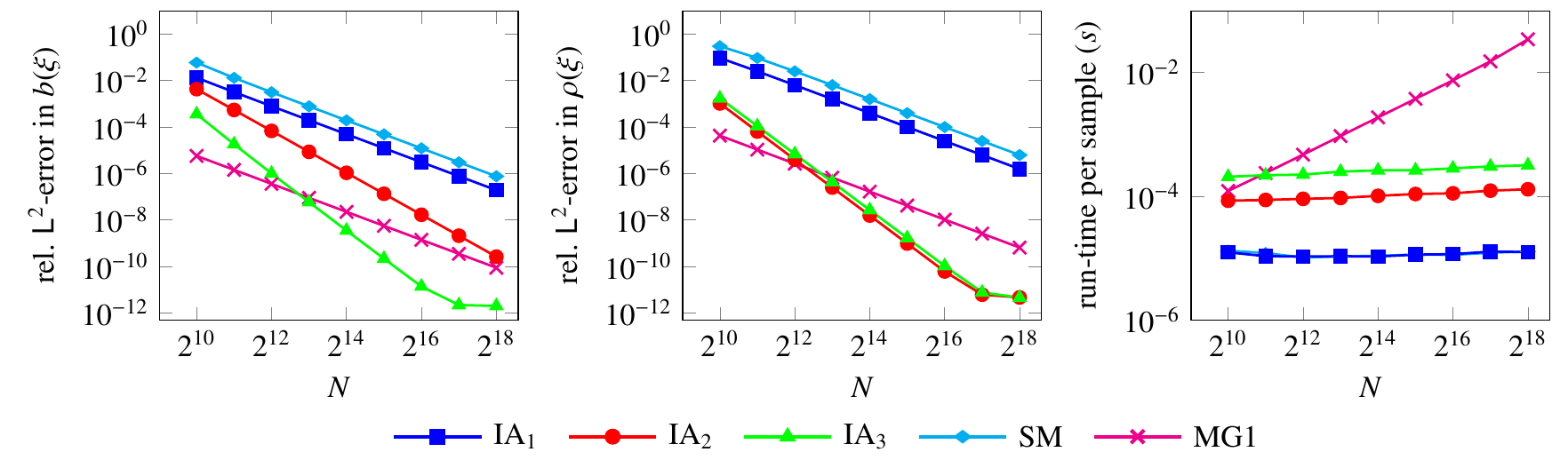}
\caption{\label{fig:results-IA} The figure shows a comparison of convergence
behavior and run-time of NFT algorithms based on the discretization schemes, namely, 
IA$_m$ ($m\in\{1,2,3\}$), Split-Magnus (SM)
and Magnus method with one-point Gauss quadrature (MG1) (the latter two 
are discussed in~\cite[Sec.~IV]{V2017INFT1} as a way of benchmarking).
The method IA$_1$ is identical to the trapezoidal rule (TR). The test
corresponds to a secant-hyperbolic profile $q(t)=4.4\sech(t)$.}
\end{figure*}
In order to develop the numerical scheme based on the implicit Adams (IA) method, 
we begin with the transformation $\tilde{\vv{v}}=e^{i\sigma_3\zeta t}\vv{v}$
so that the ZS problem in~\eqref{eq:zs-prob} reads as
\begin{equation}\label{eq:exp-int}
\tilde{\vv{v}}_t=\wtilde{U}\tilde{\vv{v}},\quad
\wtilde{U}
=\begin{pmatrix}
0 & qe^{2i\zeta t}\\
re^{-2i\zeta t} & 0
\end{pmatrix}.
\end{equation}
Let the grid $\{t_n\}$ be as defined in~\ref{sec:discrete-system}, and, set
$U_{n}=U(t_n)$ and $\wtilde{U}_{n}=\wtilde{U}(t_n)$. The discretization 
of~\eqref{eq:exp-int} using the
$m$-step IA method ($m\in\{1,2,3\}$) reads as
\begin{equation}
\tilde{\vv{v}}_{n+m}-\tilde{\vv{v}}_{n+m-1}=
h\sum_{s=0}^m\beta_s\wtilde{U}_{n+s}\tilde{\vv{v}}_{n+s},
\end{equation}
where $\vs{\beta} = (\beta_0,\beta_1,\ldots,\beta_m)$ are known 
constants~\cite[Chap.~III.1]{HNW1993} 
(also summarized in Table~\ref{tb:IA}). Solving for $\tilde{\vv{v}}_{n+m}$, we have
\begin{multline*}
\tilde{\vv{v}}_{n+m}
=\left(\sigma_0 - h\beta_m\wtilde{U}_{n+m}\right)^{-1}\times\\
\biggl[
\left(\sigma_0+h\beta_{m-1}\wtilde{U}_{n+m-1}\right)\tilde{\vv{v}}_{n+m-1}\\
+\sum_{s=0}^{m-2}h\beta_s\wtilde{U}_{n+s}\tilde{\vv{v}}_{n+s}\biggl],
\end{multline*}
or, equivalently,
\begin{multline}
{\vv{v}}_{n+m}=\left(\sigma_0 - h\beta_m{U}_{n+m}\right)^{-1}\times\\
\biggl[\sum_{s=0}^{m-2}h\beta_se^{-i\sigma_3\zeta h(m-s)}{U}_{n+s}{\vv{v}}_{n+s}\\
+e^{-i\sigma_3\zeta h}
\left(\sigma_0+h\beta_{m-1}{U}_{n+m-1}\right){\vv{v}}_{n+m-1}\biggl],
\end{multline}
where $\sigma_0=\diag(1,1)$. The individual matrices can be worked out as 
\begin{equation}
\begin{split}
&\left(\sigma_0-h\beta_m{U}_{n+m}\right)^{-1}e^{-i\sigma_3\zeta h}
\left(\sigma_0+h\beta_{m-1}{U}_{n+m-1}\right)\\
&=\frac{z^{-1}}{\Theta_{n+m}}\times\\
&\begin{pmatrix}
1+z^{2}\bar{\beta}_{m-1}R_{n+m-1}Q_{n+m} &z^{2}Q_{n+m} + \bar{\beta}_{m-1}Q_{n+m-1}\\
R_{n+m}+z^{2}\bar{\beta}_{m-1}R_{n+m-1} &
z^{2}+\bar{\beta}_{m-1}R_{n+m}Q_{n+m-1}
\end{pmatrix}\\
&\equiv z^{-1}M^{(1)}_{n+m}(z^2),
\end{split}
\end{equation}
where $Q_{n}=(h\beta_m)q_{n}$, $R_{n}=(h\beta_m)r_{n}$, $\Theta_n =
1-Q_{n}R_{n}$ and
\begin{equation}
\ovl{\vs{\beta}} = \vs{\beta}/\beta_m=(\ovl{\beta}_0,\ovl{\beta}_1,\ldots,1). 
\end{equation}
Also,
\begin{equation}
\begin{split}
&\left(\sigma_0-h\beta_m{U}_{n+m}\right)^{-1}e^{-i\sigma_3\zeta(m-s)h}h\beta_sU_{n+s}\\
&=\ovl{\beta}_s\frac{z^{-(m-s)}}{\Theta_{n+m}}\begin{pmatrix}
    z^{2(m-s)}R_{n+s}Q_{n+m} &Q_{n+s}\\
    z^{2(m-s)}R_{n+s} & R_{n+m}Q_{n+s}
\end{pmatrix}\\
&\equiv \ovl{\beta}_sz^{-(m-s)}M^{(m-s)}_{n+m}(z^2).
\end{split}
\end{equation}
\begin{table}
\def\arraystretch{1.75}
\caption{Implicit Adams Method\label{tb:IA}}
\begin{center}
\begin{tabular}{ccc}
Method & $\vs{\beta}$ & Order of Convergence\\
\hline
IA$_1$ & $(\frac{1}{2},\frac{1}{2})$ & $2$\\ 
\hline
IA$_2$ & $(-\frac{1}{12},\frac{8}{12},\frac{5}{12})$ & $3$\\ 
\hline
IA$_3$ & $(\frac{1}{24},-\frac{5}{24},\frac{19}{24},\frac{9}{24})$ &$4$\\ 
\hline\\
\end{tabular}
\end{center}
\end{table}

The $m$-step IA methods lead to transfer matrices
$\mathcal{M}_{n}\in\field{C}^{2m\times 2m}$ of the form 
\begin{equation}
\begin{split}
&{\mathcal{M}}_{n+m}(z^2) = \\
&\begin{pmatrix}
M^{(1)}_{n+m}&\ovl{\beta}_{m-2}M^{(2)}_{n+m}&
\ldots&\ovl{\beta}_{1}M^{(m-1)}_{n+m}(z^2)&\ovl{\beta}_{0}M^{(m)}_{n+m}\\
\sigma_0 &0&\ldots&0&0\\
0&\sigma_0&\ldots&0&0\\
\vdots&\vdots &\ddots & \vdots& \vdots\\
0&0&\ldots&\sigma_0&0
\end{pmatrix},
\end{split}
\end{equation}
where $M^{(s)}_{n+m}(z^2)\in\field{C}^{2\times2}$ so that
\begin{equation}
\pmb{\mathcal{W}}_{n+m}={\mathcal{M}}_{n+m}(z^2)\pmb{\mathcal{W}}_{n+m-1},
\end{equation}
where $\vv{w}_n=z^{n}\vv{v}_n$ and 
\[
\pmb{\mathcal{W}}_{n}=(\vv{w}_{n},\vv{w}_{n-1},\ldots,\vv{w}_{n-m+1})^{\tp}\in\field{C}^{2m}.
\]

Let us consider the Jost solution $\vs{\phi}(t;\zeta)$. We assume that $q_n=0$ for 
$n=-m+1, -m+2,\ldots,0$ so that 
$\vs{\phi}_n=z^{\ell_-}z^{-n}(1,0)^{\tp}$ for $n=-m+1,-m+2,\ldots,0$. The discrete 
approximation to the Jost solution can be expressed as 
$\vs{\phi}_n = z^{\ell_-}z^{-n}\vv{P}_n(z^2)$. The initial condition
works out to be 
\begin{equation*}
\begin{split}
\pmb{\mathcal{W}}_{0}
&=z^{\ell_-}
\begin{pmatrix}
\vs{\phi}_{0}\\
z\vs{\phi}_{-1}\\
\vdots\\
z^{-m+1}\vs{\phi}_{-m+1}
\end{pmatrix}\\
&=z^{\ell_-}
\begin{pmatrix}
\vs{P}_{0}(z^2)\\
\vs{P}_{-1}(z^2)\\
\vdots\\
\vs{P}_{-m+1}(z^2)
\end{pmatrix}
=z^{\ell_-}
\begin{pmatrix}
1\\
0\\
\vdots\\
1\\
0
\end{pmatrix},
\end{split}
\end{equation*}
yielding the recurrence relation
\begin{equation}
\pmb{\mathcal{P}}_{n+m}(z^2)={\mathcal{M}}_{n+m}(z^2)\pmb{\mathcal{P}}_{n+m-1}(z^2),
\end{equation}
where 
$\pmb{\mathcal{P}}_{n}(z^2) =
(\vv{P}_{n}(z^2),\vv{P}_{n-1}(z^2),\ldots,\vv{P}_{n-m+1}(z^2))^{\tp}\in\field{C}^{2m}$.
The discrete approximation to the scattering coefficients is obtained from the scattered
field: $\vs{\phi}_{N}=(a_{N} z^{-\ell_+},b_{N} z^{\ell_+})^{\tp}$ yields
$a_{N}(z^2)={P}^{(N)}_1(z^2)$ and
$b_{N}(z^2)=(z^2)^{-\ell_{+}}{P}^{(N)}_2(z^2)$. The quantities $a_{N}$ and $b_{N}$ are referred to as the 
\emph{discrete scattering coefficients} uniquely defined for $\Re\zeta\in
[-{\pi}/{2h},\,{\pi}/{2h}]$.

Finally, let us mention that, for $\zeta$ varying over a compact domain, the
error in the computation of the scattering coefficients can be shown to be
$\bigO{N^{-p}}$ provided that $q(t)$ is at least $p$-times differentiable~\cite[Chap.~III]{HNW1993}.

It is evident from the preceding paragraph that
the forward scattering step requires forming the cumulative product:
${\mathcal{M}}_{N}(z^2)\times{\mathcal{M}}_{N-1}(z^2)\times\ldots
\times{\mathcal{M}}_{2}(z^2)\times{\mathcal{M}}_{1}(z^2)$. Let $\bar{m}$ 
denote the nearest base-$2$ number greater than or equal
to $(m +1)$, then pairwise multiplication using FFT~\cite{Henrici1993} yields the 
recurrence relation for the complexity $\varpi(n)$ of computing the scattering
coefficients with $n$ samples: 
$\varpi(n) = 8m^3\nu(\bar{m}n/2)+2\varpi(n/2),\,\,n=2,\,4,\,\ldots,\,N,$ 
where $\nu(n)=\bigO{n\log n}$ is the cost of multiplying two polynomials of
degree $n-1$ (ignoring the cost of additions). Solving the
recurrence relation yields $\varpi(N)=\bigO{m^3N\log^2N}$.

Finally, the results of the tests for benchmarking are shown
Fig.~\ref{fig:results-IA}.

\section{An extension of the theorem of Epstein}
\label{app:epstein}
In the following, we would like to extend Theorem~$4$ of~\cite{Epstein2004} to
obtain the result~\eqref{eq:epstien}. Define 
the nonlinear impulse response
\begin{equation}
p(\tau) = \fourier^{-1}[\rho](\tau)
=\frac{1}{2\pi}\int_{-\infty}^{\infty}\rho(\xi) e^{-i\xi\tau}d\xi,
\end{equation}
and assume $p(\tau)\in\fs{L}^{1}\cap\fs{L}^{2}$. Consider the
Jost solutions with prescribed asymptotic behavior as 
$x\rightarrow\infty$:
\begin{equation}
\vs{\psi}(t;\zeta)=
\begin{pmatrix}
0\\
1
\end{pmatrix}e^{i\zeta t}
+\int_t^{\infty}e^{i\zeta s}{\vv{A}}(t,s)ds,
\end{equation}
where $\vv{A}$ is independent of $\zeta$. Our starting point for the analysis of
the inverse problem would be the Gelfand-Levitan-Marchenko (GLM) integral
equations. In the following we fix
$t\in\field{R}$ so that the GLM equations for $y\in\Omega_t=[t,\infty)$ is given by 
\begin{equation}\label{eq;GLM-start}
\begin{split}
&{A}_2^*(t,y)=-\int_{t}^{\infty}{A}_1(t,s){f}(s+y)ds,\\
&{A}_1^*(t,y) = f(t+y) +\int_{t}^{\infty}A_2(t,s){f}(s+y)ds,
\end{split}
\end{equation}
where $f(\tau)=p(-\tau)$. The solution of the GLM equations allows us to recover
the scattering potential using $q(t)=-2A_1(t,t)$ together with the estimate 
$\|q\chi_{[t,\infty)}\|^2_2=-2A_2(t,t)$ where $\chi_{\Omega}$ denotes the
characteristic function of $\Omega\subset\field{R}$. Define the operator
\begin{equation}\label{eq:p-GLM-OP}
\OP{P}[g](y)=\int^{\infty}_t f(y+s)g(s)ds,
\end{equation}
whose Hermitian conjugate, denoted by $\OP{P}^{\dagger}$, works out to be
\begin{equation}
\OP{P}^{\dagger}[g](y)=\int_t^{\infty}f^*(y+s)g(s)ds.
\end{equation}
Define $\OP{K}=\OP{P}^{\dagger}\circ\OP{P}$, so that
\begin{equation}
\begin{split}
\OP{K}[g](y)
&=\int_t^{\infty}ds\int_t^{\infty}dx\,f^*(y+s)f(s+x)g(x)\\
&=\int_t^{\infty}\mathcal{K}(y,x;t)g(x)dx,
\end{split}
\end{equation}
where the kernel function $\mathcal{K}(y,x;t)$ is given by
\begin{equation}
\mathcal{K}(y,x;t)
=\int_t^{\infty}ds\,f^*(y+s)f(s+x).
\end{equation}
The GLM equations in~\eqref{eq;GLM-start} can now be stated as
\begin{equation}\label{eq:fredholm}
{A}_j(t,y)={\Phi}_j(t,y)-\OP{K}[{A}_j(t,\cdot)](y),\quad j=1,2,
\end{equation}
which is a Fredholm integral equation of the second kind where
\begin{equation}
{\Phi}_1(t,y)=f^*(t+y),\quad
{\Phi}_2(t,y)=-\OP{P}^{\dagger}[f(t+\cdot)](y).
\end{equation}
Let $\mathcal{I}_m(t)=\|f\chi_{[2t,\infty)}\|_{\fs{L}^m}$ for $m=1,2,\infty$, then
\begin{equation}
\begin{split}
\|\OP{K}\|_{\fs{L}^{\infty}(\Omega_t)}
&=\esssup_{y\in\Omega_t}\int_{t}^{\infty}dx\,|\mathcal{K}(y,x;t)|\\
&\leq\esssup_{y\in\Omega_t}\int_{t}^{\infty}dx\,\int_t^{\infty}ds\,|f(y+s)||f(s+x)|\\
&\leq\esssup_{y\in\Omega_t}\int_{t+y}^{\infty}du|f(u)|\,\int_{t+u-y}^{\infty}du_1\,|f(u_1)|\\
&\leq [\mathcal{I}_1(t)]^2,
\end{split}
\end{equation}
and, $\|\Phi_2(t,\cdot)\|_{\fs{L}^{\infty}(\Omega_t)}\leq[\mathcal{I}_{2}(t)]^2$.
If $\mathcal{I}_1(t)<1$, then the standard theory of Fredholm equations suggests that the
resolvent of the operator $\OP{K}$ exists~\cite{GLS1990}. Under this assumption,
certain estimates for $q(t)$ can be easily obtained~\cite{Epstein2004}: From~\eqref{eq:fredholm}, 
we have
\begin{multline*}
\|{A}_j(t,\cdot)\|_{\fs{L}^{\infty}(\Omega_t)}\leq
\|\Phi_j(t,\cdot)\|_{\fs{L}^{\infty}(\Omega_t)}\\
+\|\OP{K}\|_{\fs{L}^{\infty}(\Omega_t)}\|{A}_j(t,\cdot)\|_{\fs{L}^2(\Omega_t)},
\end{multline*}
which yields
\begin{equation*}
\begin{split}
&\|{A}_1(t,\cdot)\|_{\fs{L}^{\infty}(\Omega_t)}\leq\frac{\mathcal{I}_{\infty}(t)}{[1-\mathcal{I}^2_1(t)]},\\
&\|{A}_2(t,\cdot)\|_{\fs{L}^{\infty}(\Omega_t)}\leq\frac{\mathcal{I}^2_2(t)}{[1-\mathcal{I}^2_1(t)]}.
\end{split}
\end{equation*}
Given that from here one can only assert that 
$|A_j(t,y)|\leq\|{A}_j(t,\cdot)\|_{\fs{L}^{\infty}(\Omega_t)}$ almost everywhere
(a.e.), we need to ascertain the continuity of $A_j(t,y)$ with respect to $y$ 
throughout the domain $\Omega_t$
or as $y\rightarrow t$ from above. Assume that $f(\tau)$ is continuous, then
$\Phi_j(t,y)$ is continuous with respect to $y$. It can be seen that the kernel function
$\mathcal{K}(y,x;t)$ is also continuous with respect to $y$. Therefore, if the resolvent 
kernel is continuous (w.r.t. $y$) then the result 
follows. To this end, consider the Neumann series for the resolvent
$\OP{R}=\sum_{n\in\field{Z}_+}(-1)^n\OP{K}_n$ where
$\OP{K}_n=\OP{K}\circ\OP{K}_{n-1}$ with $\OP{K}_1=\OP{K}$. For fixed
$t$, the partial sums 
$\sum_{1\leq n\leq N}\|\OP{K}_n\|_{\fs{L}^{\infty}(\Omega_t)}
\leq[1-\mathcal{I}^2_1(t)]^{-1}$ for all $N<\infty$.
Therefore, uniform convergence of the partial sums allows us to conclude the continuity of the 
limit of the partial sums.

Now using the identities $q(t)=-2A_1(t,t)$ and $\|q\chi_{[t,\infty)}\|^2_2=-2A_2(t,t)$, we have
\begin{equation}\label{eq:epstein-estimates}
\begin{split}
&\|q\chi_{[t,\infty)}\|_{\fs{L}^{\infty}}
\leq\frac{2\mathcal{I}_{\infty}(t)}{[1-\mathcal{I}^2_1(t)]},\\
&\|q\chi_{[t,\infty)}\|^2_{\fs{L}^{2}}
\leq\frac{2\mathcal{I}^2_2(t)}{[1-\mathcal{I}^2_1(t)]}.
\end{split}
\end{equation}
If $\mathcal{I}_1(t)<1$ does not hold for all $t\in\field{R}$, one can find a $T>0$ such 
that $\mathcal{I}_1(t)<1$ holds for $t\in[T,\infty)$. The estimates obtained
above would then be valid in $[T,\infty)$.

The second inequality in~\eqref{eq:epstein-estimates} can be used to choose the 
computational domain for the inverse
NFT. Let us consider the example considered in Sec.~\ref{sec:nbs-test}: The nonlinear impulse
response works out to be
\begin{equation}
p_{\text{rc}}(\tau) = \frac{A}{\pi\tau_s}\sinc\left(\frac{\tau}{\tau_s}\right)
\frac{\cos\left(\beta\frac{\tau}{\tau_s}\right)}{1-\left(\frac{2\beta\tau}{\pi\tau_s}\right)^2}.
\end{equation}
Note that $p_{\text{rc}}(-\tau)=p_{\text{rc}}(\tau)$. From the asymptotic form 
\[
|p_{\text{rc}}(\tau)|\sim\left(\frac{A\pi\tau_s^2}{4\beta^2}\right)\frac{1}{\tau^3},
\]
it follows that 
\begin{equation}
\begin{split}
&\mathcal{I}_1(T)\sim\left(\frac{A\pi\tau_s^2}{4\beta^2}\right)\frac{1}{2(2T)^2},\\
&\mathcal{I}^2_2(T)\sim\left(\frac{A\pi\tau_s^2}{4\beta^2}\right)^2\frac{1}{5(2T)^5}.
\end{split}
\end{equation}
If $\mathcal{I}_1(T)\ll 1$, then setting $\mathcal{I}^2_2(T)=\epsilon$ gives
\begin{equation}
T(\epsilon)\sim
\frac{1}{2}\left(\frac{\pi^2A^2\tau_s^4}{40\beta^4\epsilon}\right)^{1/5}.
\end{equation}

\end{document}